\documentclass[a4paper]{article}
\usepackage[utf8]{inputenc}
\usepackage[T1]{fontenc} % ou \usepackage[OT1]{fontenc}
\usepackage{RR,RRthemes}
\usepackage[%
  bookmarks, 
  linkcolor= blue,
  citecolor= blue,
  filecolor= blue,
  urlcolor=  blue,
  colorlinks=true, 
  hyperindex=true,
  pdftitle={A Markov model of land use dynamics}, 
  pdfauthor={Fabien Campillo, Dominique Herv\'e, Angelo Raherinirina, and Rivo Rakotozafy}
]{hyperref}             % Producing hyperlinked PDF
\usepackage[english]{babel}
\usepackage{graphics}
\usepackage{graphicx}           
\usepackage{url}
\usepackage{makeidx}
\usepackage[plain]{algorithm}
\usepackage{algorithmic} 
\usepackage{color} 
\usepackage{rotating}
\usepackage[table]{xcolor}
\usepackage[margin=3.3cm]{geometry}
\usepackage{pifont}
\usepackage{tikz}
\usepackage{amsmath,amsfonts,amssymb,latexsym,stmaryrd}
\usepackage{multirow}
\newtheorem{theorem*}     {theorem}

\newcommand{\eqdef}     {\stackrel{\textup{\tiny def}}{=}}

\newsavebox{\fmbox}

\newcommand{\var}{{\textrm{var}}}

\newcommand{\F}        {\mathbb F}

\newcommand{\N}        {\mathbb N}

\newcommand{\E}        {\mathbb E}

\renewcommand{\P}      {\mathbb P}

\newcommand{\simiid}     {\stackrel{\textrm{iid}}{\sim}}

\newcommand{\ot}        {\leftarrow}

\newcommand{\II}   {{\mathcal I}}

\newcommand{\NN}   {{\mathcal N}}

\newcommand{\rmd}   {{{\textrm{\upshape d}}}}

\newcommand{\indic}{{\mathrm\mathbf1}}
\renewcommand{\epsilon}{\varepsilon}

\def\dobm{
    \copy1\kern-\wd1\kern0.05ex\copy1\kern-\wd1\kern0.05ex\box1}

\newcommand{\fenumi}  {\textrm{\rm({\textit{i}}\/)}}
\newcommand{\fenumii} {\textrm{\rm({\textit{ii}}\/)}}
\newcommand{\fenumiii}{\textrm{\rm({\textit{iii}}\/)}}
\newcommand{\fenumiv} {\textrm{\rm({\textit{iv}}\/)}}

\newcommand
      {\sysdys}
      {{\sf S\kern-.15em\raise.3ex\hbox{Y}\kern-.15em
            SD\kern-.15em\raise.3ex\hbox{Y}\kern-.15emS}}

\colorlet{darkcyan}  {cyan!80!black}
\colorlet{darkred}  {red!80!black}
\colorlet{darkblue} {blue!80!black}
\colorlet{darkgreen}{green!50!black}

\definecolor{Fsnow}{rgb}{1.000,0.980,0.980}
\definecolor{FGhostWhite}{rgb}{0.973,0.973,1.000}
\definecolor{FWhiteSmoke}{rgb}{0.961,0.961,0.961}
\definecolor{Fgainsboro}{rgb}{0.863,0.863,0.863}
\definecolor{FFloralWhite}{rgb}{1.000,0.980,0.941}
\definecolor{FOldLace}{rgb}{0.992,0.961,0.902}
\definecolor{Flinen}{rgb}{0.980,0.941,0.902}
\definecolor{FAntiqueWhite}{rgb}{0.980,0.922,0.843}
\definecolor{FPapayaWhip}{rgb}{1.000,0.937,0.835}
\definecolor{FBlanchedAlmond}{rgb}{1.000,0.922,0.804}
\definecolor{Fbisque}{rgb}{1.000,0.894,0.769}
\definecolor{FPeachPuff}{rgb}{1.000,0.855,0.725}
\definecolor{FNavajoWhite}{rgb}{1.000,0.871,0.678}
\definecolor{Fmoccasin}{rgb}{1.000,0.894,0.710}
\definecolor{Fcornsilk}{rgb}{1.000,0.973,0.863}
\definecolor{Fivory}{rgb}{1.000,1.000,0.941}
\definecolor{FLemonChiffon}{rgb}{1.000,0.980,0.804}
\definecolor{Fseashell}{rgb}{1.000,0.961,0.933}
\definecolor{Fhoneydew}{rgb}{0.941,1.000,0.941}
\definecolor{FMintCream}{rgb}{0.961,1.000,0.980}
\definecolor{Fazure}{rgb}{0.941,1.000,1.000}
\definecolor{FAliceBlue}{rgb}{0.941,0.973,1.000}
\definecolor{Flavender}{rgb}{0.902,0.902,0.980}
\definecolor{FLavenderBlush}{rgb}{1.000,0.941,0.961}
\definecolor{FMistyRose}{rgb}{1.000,0.894,0.882}
\definecolor{Fwhite}{rgb}{1.000,1.000,1.000}
\definecolor{Fblack}{rgb}{0.000,0.000,0.000}
\definecolor{FDarkSlateGray}{rgb}{0.184,0.310,0.310}
\definecolor{FDarkSlateGrey}{rgb}{0.184,0.310,0.310}
\definecolor{FDimGray}{rgb}{0.412,0.412,0.412}
\definecolor{FDimGrey}{rgb}{0.412,0.412,0.412}
\definecolor{FSlateGray}{rgb}{0.439,0.502,0.565}
\definecolor{FSlateGrey}{rgb}{0.439,0.502,0.565}
\definecolor{FLightSlateGray}{rgb}{0.467,0.533,0.600}
\definecolor{FLightSlateGrey}{rgb}{0.467,0.533,0.600}
\definecolor{Fgray}{rgb}{0.745,0.745,0.745}
\definecolor{Fgrey}{rgb}{0.745,0.745,0.745}
\definecolor{FLightGrey}{rgb}{0.827,0.827,0.827}
\definecolor{FLightGray}{rgb}{0.827,0.827,0.827}
\definecolor{FMidnightBlue}{rgb}{0.098,0.098,0.439}
\definecolor{Fnavy}{rgb}{0.000,0.000,0.502}
\definecolor{FNavyBlue}{rgb}{0.000,0.000,0.502}
\definecolor{FCornflowerBlue}{rgb}{0.392,0.584,0.929}
\definecolor{FDarkSlateBlue}{rgb}{0.282,0.239,0.545}
\definecolor{FSlateBlue}{rgb}{0.416,0.353,0.804}
\definecolor{FMediumSlateBlue}{rgb}{0.482,0.408,0.933}
\definecolor{FLightSlateBlue}{rgb}{0.518,0.439,1.000}
\definecolor{FMediumBlue}{rgb}{0.000,0.000,0.804}
\definecolor{FRoyalBlue}{rgb}{0.255,0.412,0.882}
\definecolor{Fblue}{rgb}{0.000,0.000,1.000}
\definecolor{FDodgerBlue}{rgb}{0.118,0.565,1.000}
\definecolor{FDeepSkyBlue}{rgb}{0.000,0.749,1.000}
\definecolor{FSkyBlue}{rgb}{0.529,0.808,0.922}
\definecolor{FLightSkyBlue}{rgb}{0.529,0.808,0.980}
\definecolor{FSteelBlue}{rgb}{0.275,0.510,0.706}
\definecolor{FLightSteelBlue}{rgb}{0.690,0.769,0.871}
\definecolor{FLightBlue}{rgb}{0.678,0.847,0.902}
\definecolor{FPowderBlue}{rgb}{0.690,0.878,0.902}
\definecolor{FPaleTurquoise}{rgb}{0.686,0.933,0.933}
\definecolor{FDarkTurquoise}{rgb}{0.000,0.808,0.820}
\definecolor{FMediumTurquoise}{rgb}{0.282,0.820,0.800}
\definecolor{Fturquoise}{rgb}{0.251,0.878,0.816}
\definecolor{Fcyan}{rgb}{0.000,1.000,1.000}
\definecolor{FLightCyan}{rgb}{0.878,1.000,1.000}
\definecolor{FCadetBlue}{rgb}{0.373,0.620,0.627}
\definecolor{FMediumAquamarine}{rgb}{0.400,0.804,0.667}
\definecolor{Faquamarine}{rgb}{0.498,1.000,0.831}
\definecolor{FDarkGreen}{rgb}{0.000,0.392,0.000}
\definecolor{FDarkOliveGreen}{rgb}{0.333,0.420,0.184}
\definecolor{FDarkSeaGreen}{rgb}{0.561,0.737,0.561}
\definecolor{FSeaGreen}{rgb}{0.180,0.545,0.341}
\definecolor{FMediumSeaGreen}{rgb}{0.235,0.702,0.443}
\definecolor{FLightSeaGreen}{rgb}{0.125,0.698,0.667}
\definecolor{FPaleGreen}{rgb}{0.596,0.984,0.596}
\definecolor{FSpringGreen}{rgb}{0.000,1.000,0.498}
\definecolor{FLawnGreen}{rgb}{0.486,0.988,0.000}
\definecolor{Fgreen}{rgb}{0.000,1.000,0.000}
\definecolor{Fchartreuse}{rgb}{0.498,1.000,0.000}
\definecolor{FMediumSpringGreen}{rgb}{0.000,0.980,0.604}
\definecolor{FGreenYellow}{rgb}{0.678,1.000,0.184}
\definecolor{FLimeGreen}{rgb}{0.196,0.804,0.196}
\definecolor{FYellowGreen}{rgb}{0.604,0.804,0.196}
\definecolor{FForestGreen}{rgb}{0.133,0.545,0.133}
\definecolor{FOliveDrab}{rgb}{0.420,0.557,0.137}
\definecolor{FDarkKhaki}{rgb}{0.741,0.718,0.420}
\definecolor{Fkhaki}{rgb}{0.941,0.902,0.549}
\definecolor{FPaleGoldenrod}{rgb}{0.933,0.910,0.667}
\definecolor{FLightGoldenrodYellow}{rgb}{0.980,0.980,0.824}
\definecolor{FLightYellow}{rgb}{1.000,1.000,0.878}
\definecolor{Fyellow}{rgb}{1.000,1.000,0.000}
\definecolor{Fgold}{rgb}{1.000,0.843,0.000}
\definecolor{FLightGoldenrod}{rgb}{0.933,0.867,0.510}
\definecolor{Fgoldenrod}{rgb}{0.855,0.647,0.125}
\definecolor{FDarkGoldenrod}{rgb}{0.722,0.525,0.043}
\definecolor{FRosyBrown}{rgb}{0.737,0.561,0.561}
\definecolor{FIndianRed}{rgb}{0.804,0.361,0.361}
\definecolor{FSaddleBrown}{rgb}{0.545,0.271,0.075}
\definecolor{Fsienna}{rgb}{0.627,0.322,0.176}
\definecolor{Fperu}{rgb}{0.804,0.522,0.247}
\definecolor{Fburlywood}{rgb}{0.871,0.722,0.529}
\definecolor{Fbeige}{rgb}{0.961,0.961,0.863}
\definecolor{Fwheat}{rgb}{0.961,0.871,0.702}
\definecolor{FSandyBrown}{rgb}{0.957,0.643,0.376}
\definecolor{Ftan}{rgb}{0.824,0.706,0.549}
\definecolor{Fchocolate}{rgb}{0.824,0.412,0.118}
\definecolor{Ffirebrick}{rgb}{0.698,0.133,0.133}
\definecolor{Fbrown}{rgb}{0.647,0.165,0.165}
\definecolor{FDarkSalmon}{rgb}{0.914,0.588,0.478}
\definecolor{Fsalmon}{rgb}{0.980,0.502,0.447}
\definecolor{FLightSalmon}{rgb}{1.000,0.627,0.478}
\definecolor{Forange}{rgb}{1.000,0.647,0.000}
\definecolor{FDarkOrange}{rgb}{1.000,0.549,0.000}
\definecolor{Fcoral}{rgb}{1.000,0.498,0.314}
\definecolor{FLightCoral}{rgb}{0.941,0.502,0.502}
\definecolor{Ftomato}{rgb}{1.000,0.388,0.278}
\definecolor{FOrangeRed}{rgb}{1.000,0.271,0.000}
\definecolor{Fred}{rgb}{1.000,0.000,0.000}
\definecolor{FHotPink}{rgb}{1.000,0.412,0.706}
\definecolor{FDeepPink}{rgb}{1.000,0.078,0.576}
\definecolor{Fpink}{rgb}{1.000,0.753,0.796}
\definecolor{FLightPink}{rgb}{1.000,0.714,0.757}
\definecolor{FPaleVioletRed}{rgb}{0.859,0.439,0.576}
\definecolor{Fmaroon}{rgb}{0.690,0.188,0.376}
\definecolor{FMediumVioletRed}{rgb}{0.780,0.082,0.522}
\definecolor{FVioletRed}{rgb}{0.816,0.125,0.565}
\definecolor{Fmagenta}{rgb}{1.000,0.000,1.000}
\definecolor{Fviolet}{rgb}{0.933,0.510,0.933}
\definecolor{Fplum}{rgb}{0.867,0.627,0.867}
\definecolor{Forchid}{rgb}{0.855,0.439,0.839}
\definecolor{FMediumOrchid}{rgb}{0.729,0.333,0.827}
\definecolor{FDarkOrchid}{rgb}{0.600,0.196,0.800}
\definecolor{FDarkViolet}{rgb}{0.580,0.000,0.827}
\definecolor{FBlueViolet}{rgb}{0.541,0.169,0.886}
\definecolor{Fpurple}{rgb}{0.627,0.125,0.941}
\definecolor{FMediumPurple}{rgb}{0.576,0.439,0.859}
\definecolor{Fthistle}{rgb}{0.847,0.749,0.847}
\definecolor{Fsnow1}{rgb}{1.000,0.980,0.980}
\definecolor{Fsnow2}{rgb}{0.933,0.914,0.914}
\definecolor{Fsnow3}{rgb}{0.804,0.788,0.788}
\definecolor{Fsnow4}{rgb}{0.545,0.537,0.537}
\definecolor{Fseashell1}{rgb}{1.000,0.961,0.933}
\definecolor{Fseashell2}{rgb}{0.933,0.898,0.871}
\definecolor{Fseashell3}{rgb}{0.804,0.773,0.749}
\definecolor{Fseashell4}{rgb}{0.545,0.525,0.510}
\definecolor{FAntiqueWhite1}{rgb}{1.000,0.937,0.859}
\definecolor{FAntiqueWhite2}{rgb}{0.933,0.875,0.800}
\definecolor{FAntiqueWhite3}{rgb}{0.804,0.753,0.690}
\definecolor{FAntiqueWhite4}{rgb}{0.545,0.514,0.471}
\definecolor{Fbisque1}{rgb}{1.000,0.894,0.769}
\definecolor{Fbisque2}{rgb}{0.933,0.835,0.718}
\definecolor{Fbisque3}{rgb}{0.804,0.718,0.620}
\definecolor{Fbisque4}{rgb}{0.545,0.490,0.420}
\definecolor{FPeachPuff1}{rgb}{1.000,0.855,0.725}
\definecolor{FPeachPuff2}{rgb}{0.933,0.796,0.678}
\definecolor{FPeachPuff3}{rgb}{0.804,0.686,0.584}
\definecolor{FPeachPuff4}{rgb}{0.545,0.467,0.396}
\definecolor{FNavajoWhite1}{rgb}{1.000,0.871,0.678}
\definecolor{FNavajoWhite2}{rgb}{0.933,0.812,0.631}
\definecolor{FNavajoWhite3}{rgb}{0.804,0.702,0.545}
\definecolor{FNavajoWhite4}{rgb}{0.545,0.475,0.369}
\definecolor{FLemonChiffon1}{rgb}{1.000,0.980,0.804}
\definecolor{FLemonChiffon2}{rgb}{0.933,0.914,0.749}
\definecolor{FLemonChiffon3}{rgb}{0.804,0.788,0.647}
\definecolor{FLemonChiffon4}{rgb}{0.545,0.537,0.439}
\definecolor{Fcornsilk1}{rgb}{1.000,0.973,0.863}
\definecolor{Fcornsilk2}{rgb}{0.933,0.910,0.804}
\definecolor{Fcornsilk3}{rgb}{0.804,0.784,0.694}
\definecolor{Fcornsilk4}{rgb}{0.545,0.533,0.471}
\definecolor{Fivory1}{rgb}{1.000,1.000,0.941}
\definecolor{Fivory2}{rgb}{0.933,0.933,0.878}
\definecolor{Fivory3}{rgb}{0.804,0.804,0.757}
\definecolor{Fivory4}{rgb}{0.545,0.545,0.514}
\definecolor{Fhoneydew1}{rgb}{0.941,1.000,0.941}
\definecolor{Fhoneydew2}{rgb}{0.878,0.933,0.878}
\definecolor{Fhoneydew3}{rgb}{0.757,0.804,0.757}
\definecolor{Fhoneydew4}{rgb}{0.514,0.545,0.514}
\definecolor{FLavenderBlush1}{rgb}{1.000,0.941,0.961}
\definecolor{FLavenderBlush2}{rgb}{0.933,0.878,0.898}
\definecolor{FLavenderBlush3}{rgb}{0.804,0.757,0.773}
\definecolor{FLavenderBlush4}{rgb}{0.545,0.514,0.525}
\definecolor{FMistyRose1}{rgb}{1.000,0.894,0.882}
\definecolor{FMistyRose2}{rgb}{0.933,0.835,0.824}
\definecolor{FMistyRose3}{rgb}{0.804,0.718,0.710}
\definecolor{FMistyRose4}{rgb}{0.545,0.490,0.482}
\definecolor{Fazure1}{rgb}{0.941,1.000,1.000}
\definecolor{Fazure2}{rgb}{0.878,0.933,0.933}
\definecolor{Fazure3}{rgb}{0.757,0.804,0.804}
\definecolor{Fazure4}{rgb}{0.514,0.545,0.545}
\definecolor{FSlateBlue1}{rgb}{0.514,0.435,1.000}
\definecolor{FSlateBlue2}{rgb}{0.478,0.404,0.933}
\definecolor{FSlateBlue3}{rgb}{0.412,0.349,0.804}
\definecolor{FSlateBlue4}{rgb}{0.278,0.235,0.545}
\definecolor{FRoyalBlue1}{rgb}{0.282,0.463,1.000}
\definecolor{FRoyalBlue2}{rgb}{0.263,0.431,0.933}
\definecolor{FRoyalBlue3}{rgb}{0.227,0.373,0.804}
\definecolor{FRoyalBlue4}{rgb}{0.153,0.251,0.545}
\definecolor{Fblue1}{rgb}{0.000,0.000,1.000}
\definecolor{Fblue2}{rgb}{0.000,0.000,0.933}
\definecolor{Fblue3}{rgb}{0.000,0.000,0.804}
\definecolor{Fblue4}{rgb}{0.000,0.000,0.545}
\definecolor{FDodgerBlue1}{rgb}{0.118,0.565,1.000}
\definecolor{FDodgerBlue2}{rgb}{0.110,0.525,0.933}
\definecolor{FDodgerBlue3}{rgb}{0.094,0.455,0.804}
\definecolor{FDodgerBlue4}{rgb}{0.063,0.306,0.545}
\definecolor{FSteelBlue1}{rgb}{0.388,0.722,1.000}
\definecolor{FSteelBlue2}{rgb}{0.361,0.675,0.933}
\definecolor{FSteelBlue3}{rgb}{0.310,0.580,0.804}
\definecolor{FSteelBlue4}{rgb}{0.212,0.392,0.545}
\definecolor{FDeepSkyBlue1}{rgb}{0.000,0.749,1.000}
\definecolor{FDeepSkyBlue2}{rgb}{0.000,0.698,0.933}
\definecolor{FDeepSkyBlue3}{rgb}{0.000,0.604,0.804}
\definecolor{FDeepSkyBlue4}{rgb}{0.000,0.408,0.545}
\definecolor{FSkyBlue1}{rgb}{0.529,0.808,1.000}
\definecolor{FSkyBlue2}{rgb}{0.494,0.753,0.933}
\definecolor{FSkyBlue3}{rgb}{0.424,0.651,0.804}
\definecolor{FSkyBlue4}{rgb}{0.290,0.439,0.545}
\definecolor{FLightSkyBlue1}{rgb}{0.690,0.886,1.000}
\definecolor{FLightSkyBlue2}{rgb}{0.643,0.827,0.933}
\definecolor{FLightSkyBlue3}{rgb}{0.553,0.714,0.804}
\definecolor{FLightSkyBlue4}{rgb}{0.376,0.482,0.545}
\definecolor{FSlateGray1}{rgb}{0.776,0.886,1.000}
\definecolor{FSlateGray2}{rgb}{0.725,0.827,0.933}
\definecolor{FSlateGray3}{rgb}{0.624,0.714,0.804}
\definecolor{FSlateGray4}{rgb}{0.424,0.482,0.545}
\definecolor{FLightSteelBlue1}{rgb}{0.792,0.882,1.000}
\definecolor{FLightSteelBlue2}{rgb}{0.737,0.824,0.933}
\definecolor{FLightSteelBlue3}{rgb}{0.635,0.710,0.804}
\definecolor{FLightSteelBlue4}{rgb}{0.431,0.482,0.545}
\definecolor{FLightBlue1}{rgb}{0.749,0.937,1.000}
\definecolor{FLightBlue2}{rgb}{0.698,0.875,0.933}
\definecolor{FLightBlue3}{rgb}{0.604,0.753,0.804}
\definecolor{FLightBlue4}{rgb}{0.408,0.514,0.545}
\definecolor{FLightCyan1}{rgb}{0.878,1.000,1.000}
\definecolor{FLightCyan2}{rgb}{0.820,0.933,0.933}
\definecolor{FLightCyan3}{rgb}{0.706,0.804,0.804}
\definecolor{FLightCyan4}{rgb}{0.478,0.545,0.545}
\definecolor{FPaleTurquoise1}{rgb}{0.733,1.000,1.000}
\definecolor{FPaleTurquoise2}{rgb}{0.682,0.933,0.933}
\definecolor{FPaleTurquoise3}{rgb}{0.588,0.804,0.804}
\definecolor{FPaleTurquoise4}{rgb}{0.400,0.545,0.545}
\definecolor{FCadetBlue1}{rgb}{0.596,0.961,1.000}
\definecolor{FCadetBlue2}{rgb}{0.557,0.898,0.933}
\definecolor{FCadetBlue3}{rgb}{0.478,0.773,0.804}
\definecolor{FCadetBlue4}{rgb}{0.325,0.525,0.545}
\definecolor{Fturquoise1}{rgb}{0.000,0.961,1.000}
\definecolor{Fturquoise2}{rgb}{0.000,0.898,0.933}
\definecolor{Fturquoise3}{rgb}{0.000,0.773,0.804}
\definecolor{Fturquoise4}{rgb}{0.000,0.525,0.545}
\definecolor{Fcyan1}{rgb}{0.000,1.000,1.000}
\definecolor{Fcyan2}{rgb}{0.000,0.933,0.933}
\definecolor{Fcyan3}{rgb}{0.000,0.804,0.804}
\definecolor{Fcyan4}{rgb}{0.000,0.545,0.545}
\definecolor{FDarkSlateGray1}{rgb}{0.592,1.000,1.000}
\definecolor{FDarkSlateGray2}{rgb}{0.553,0.933,0.933}
\definecolor{FDarkSlateGray3}{rgb}{0.475,0.804,0.804}
\definecolor{FDarkSlateGray4}{rgb}{0.322,0.545,0.545}
\definecolor{Faquamarine1}{rgb}{0.498,1.000,0.831}
\definecolor{Faquamarine2}{rgb}{0.463,0.933,0.776}
\definecolor{Faquamarine3}{rgb}{0.400,0.804,0.667}
\definecolor{Faquamarine4}{rgb}{0.271,0.545,0.455}
\definecolor{FDarkSeaGreen1}{rgb}{0.757,1.000,0.757}
\definecolor{FDarkSeaGreen2}{rgb}{0.706,0.933,0.706}
\definecolor{FDarkSeaGreen3}{rgb}{0.608,0.804,0.608}
\definecolor{FDarkSeaGreen4}{rgb}{0.412,0.545,0.412}
\definecolor{FSeaGreen1}{rgb}{0.329,1.000,0.624}
\definecolor{FSeaGreen2}{rgb}{0.306,0.933,0.580}
\definecolor{FSeaGreen3}{rgb}{0.263,0.804,0.502}
\definecolor{FSeaGreen4}{rgb}{0.180,0.545,0.341}
\definecolor{FPaleGreen1}{rgb}{0.604,1.000,0.604}
\definecolor{FPaleGreen2}{rgb}{0.565,0.933,0.565}
\definecolor{FPaleGreen3}{rgb}{0.486,0.804,0.486}
\definecolor{FPaleGreen4}{rgb}{0.329,0.545,0.329}
\definecolor{FSpringGreen1}{rgb}{0.000,1.000,0.498}
\definecolor{FSpringGreen2}{rgb}{0.000,0.933,0.463}
\definecolor{FSpringGreen3}{rgb}{0.000,0.804,0.400}
\definecolor{FSpringGreen4}{rgb}{0.000,0.545,0.271}
\definecolor{Fgreen1}{rgb}{0.000,1.000,0.000}
\definecolor{Fgreen2}{rgb}{0.000,0.933,0.000}
\definecolor{Fgreen3}{rgb}{0.000,0.804,0.000}
\definecolor{Fgreen4}{rgb}{0.000,0.545,0.000}
\definecolor{Fchartreuse1}{rgb}{0.498,1.000,0.000}
\definecolor{Fchartreuse2}{rgb}{0.463,0.933,0.000}
\definecolor{Fchartreuse3}{rgb}{0.400,0.804,0.000}
\definecolor{Fchartreuse4}{rgb}{0.271,0.545,0.000}
\definecolor{FOliveDrab1}{rgb}{0.753,1.000,0.243}
\definecolor{FOliveDrab2}{rgb}{0.702,0.933,0.227}
\definecolor{FOliveDrab3}{rgb}{0.604,0.804,0.196}
\definecolor{FOliveDrab4}{rgb}{0.412,0.545,0.133}
\definecolor{FDarkOliveGreen1}{rgb}{0.792,1.000,0.439}
\definecolor{FDarkOliveGreen2}{rgb}{0.737,0.933,0.408}
\definecolor{FDarkOliveGreen3}{rgb}{0.635,0.804,0.353}
\definecolor{FDarkOliveGreen4}{rgb}{0.431,0.545,0.239}
\definecolor{Fkhaki1}{rgb}{1.000,0.965,0.561}
\definecolor{Fkhaki2}{rgb}{0.933,0.902,0.522}
\definecolor{Fkhaki3}{rgb}{0.804,0.776,0.451}
\definecolor{Fkhaki4}{rgb}{0.545,0.525,0.306}
\definecolor{FLightGoldenrod1}{rgb}{1.000,0.925,0.545}
\definecolor{FLightGoldenrod2}{rgb}{0.933,0.863,0.510}
\definecolor{FLightGoldenrod3}{rgb}{0.804,0.745,0.439}
\definecolor{FLightGoldenrod4}{rgb}{0.545,0.506,0.298}
\definecolor{FLightYellow1}{rgb}{1.000,1.000,0.878}
\definecolor{FLightYellow2}{rgb}{0.933,0.933,0.820}
\definecolor{FLightYellow3}{rgb}{0.804,0.804,0.706}
\definecolor{FLightYellow4}{rgb}{0.545,0.545,0.478}
\definecolor{Fyellow1}{rgb}{1.000,1.000,0.000}
\definecolor{Fyellow2}{rgb}{0.933,0.933,0.000}
\definecolor{Fyellow3}{rgb}{0.804,0.804,0.000}
\definecolor{Fyellow4}{rgb}{0.545,0.545,0.000}
\definecolor{Fgold1}{rgb}{1.000,0.843,0.000}
\definecolor{Fgold2}{rgb}{0.933,0.788,0.000}
\definecolor{Fgold3}{rgb}{0.804,0.678,0.000}
\definecolor{Fgold4}{rgb}{0.545,0.459,0.000}
\definecolor{Fgoldenrod1}{rgb}{1.000,0.757,0.145}
\definecolor{Fgoldenrod2}{rgb}{0.933,0.706,0.133}
\definecolor{Fgoldenrod3}{rgb}{0.804,0.608,0.114}
\definecolor{Fgoldenrod4}{rgb}{0.545,0.412,0.078}
\definecolor{FDarkGoldenrod1}{rgb}{1.000,0.725,0.059}
\definecolor{FDarkGoldenrod2}{rgb}{0.933,0.678,0.055}
\definecolor{FDarkGoldenrod3}{rgb}{0.804,0.584,0.047}
\definecolor{FDarkGoldenrod4}{rgb}{0.545,0.396,0.031}
\definecolor{FRosyBrown1}{rgb}{1.000,0.757,0.757}
\definecolor{FRosyBrown2}{rgb}{0.933,0.706,0.706}
\definecolor{FRosyBrown3}{rgb}{0.804,0.608,0.608}
\definecolor{FRosyBrown4}{rgb}{0.545,0.412,0.412}
\definecolor{FIndianRed1}{rgb}{1.000,0.416,0.416}
\definecolor{FIndianRed2}{rgb}{0.933,0.388,0.388}
\definecolor{FIndianRed3}{rgb}{0.804,0.333,0.333}
\definecolor{FIndianRed4}{rgb}{0.545,0.227,0.227}
\definecolor{Fsienna1}{rgb}{1.000,0.510,0.278}
\definecolor{Fsienna2}{rgb}{0.933,0.475,0.259}
\definecolor{Fsienna3}{rgb}{0.804,0.408,0.224}
\definecolor{Fsienna4}{rgb}{0.545,0.278,0.149}
\definecolor{Fburlywood1}{rgb}{1.000,0.827,0.608}
\definecolor{Fburlywood2}{rgb}{0.933,0.773,0.569}
\definecolor{Fburlywood3}{rgb}{0.804,0.667,0.490}
\definecolor{Fburlywood4}{rgb}{0.545,0.451,0.333}
\definecolor{Fwheat1}{rgb}{1.000,0.906,0.729}
\definecolor{Fwheat2}{rgb}{0.933,0.847,0.682}
\definecolor{Fwheat3}{rgb}{0.804,0.729,0.588}
\definecolor{Fwheat4}{rgb}{0.545,0.494,0.400}
\definecolor{Ftan1}{rgb}{1.000,0.647,0.310}
\definecolor{Ftan2}{rgb}{0.933,0.604,0.286}
\definecolor{Ftan3}{rgb}{0.804,0.522,0.247}
\definecolor{Ftan4}{rgb}{0.545,0.353,0.169}
\definecolor{Fchocolate1}{rgb}{1.000,0.498,0.141}
\definecolor{Fchocolate2}{rgb}{0.933,0.463,0.129}
\definecolor{Fchocolate3}{rgb}{0.804,0.400,0.114}
\definecolor{Fchocolate4}{rgb}{0.545,0.271,0.075}
\definecolor{Ffirebrick1}{rgb}{1.000,0.188,0.188}
\definecolor{Ffirebrick2}{rgb}{0.933,0.173,0.173}
\definecolor{Ffirebrick3}{rgb}{0.804,0.149,0.149}
\definecolor{Ffirebrick4}{rgb}{0.545,0.102,0.102}
\definecolor{Fbrown1}{rgb}{1.000,0.251,0.251}
\definecolor{Fbrown2}{rgb}{0.933,0.231,0.231}
\definecolor{Fbrown3}{rgb}{0.804,0.200,0.200}
\definecolor{Fbrown4}{rgb}{0.545,0.137,0.137}
\definecolor{Fsalmon1}{rgb}{1.000,0.549,0.412}
\definecolor{Fsalmon2}{rgb}{0.933,0.510,0.384}
\definecolor{Fsalmon3}{rgb}{0.804,0.439,0.329}
\definecolor{Fsalmon4}{rgb}{0.545,0.298,0.224}
\definecolor{FLightSalmon1}{rgb}{1.000,0.627,0.478}
\definecolor{FLightSalmon2}{rgb}{0.933,0.584,0.447}
\definecolor{FLightSalmon3}{rgb}{0.804,0.506,0.384}
\definecolor{FLightSalmon4}{rgb}{0.545,0.341,0.259}
\definecolor{Forange1}{rgb}{1.000,0.647,0.000}
\definecolor{Forange2}{rgb}{0.933,0.604,0.000}
\definecolor{Forange3}{rgb}{0.804,0.522,0.000}
\definecolor{Forange4}{rgb}{0.545,0.353,0.000}
\definecolor{FDarkOrange1}{rgb}{1.000,0.498,0.000}
\definecolor{FDarkOrange2}{rgb}{0.933,0.463,0.000}
\definecolor{FDarkOrange3}{rgb}{0.804,0.400,0.000}
\definecolor{FDarkOrange4}{rgb}{0.545,0.271,0.000}
\definecolor{Fcoral1}{rgb}{1.000,0.447,0.337}
\definecolor{Fcoral2}{rgb}{0.933,0.416,0.314}
\definecolor{Fcoral3}{rgb}{0.804,0.357,0.271}
\definecolor{Fcoral4}{rgb}{0.545,0.243,0.184}
\definecolor{Ftomato1}{rgb}{1.000,0.388,0.278}
\definecolor{Ftomato2}{rgb}{0.933,0.361,0.259}
\definecolor{Ftomato3}{rgb}{0.804,0.310,0.224}
\definecolor{Ftomato4}{rgb}{0.545,0.212,0.149}
\definecolor{FOrangeRed1}{rgb}{1.000,0.271,0.000}
\definecolor{FOrangeRed2}{rgb}{0.933,0.251,0.000}
\definecolor{FOrangeRed3}{rgb}{0.804,0.216,0.000}
\definecolor{FOrangeRed4}{rgb}{0.545,0.145,0.000}
\definecolor{Fred1}{rgb}{1.000,0.000,0.000}
\definecolor{Fred2}{rgb}{0.933,0.000,0.000}
\definecolor{Fred3}{rgb}{0.804,0.000,0.000}
\definecolor{Fred4}{rgb}{0.545,0.000,0.000}
\definecolor{FDeepPink1}{rgb}{1.000,0.078,0.576}
\definecolor{FDeepPink2}{rgb}{0.933,0.071,0.537}
\definecolor{FDeepPink3}{rgb}{0.804,0.063,0.463}
\definecolor{FDeepPink4}{rgb}{0.545,0.039,0.314}
\definecolor{FHotPink1}{rgb}{1.000,0.431,0.706}
\definecolor{FHotPink2}{rgb}{0.933,0.416,0.655}
\definecolor{FHotPink3}{rgb}{0.804,0.376,0.565}
\definecolor{FHotPink4}{rgb}{0.545,0.227,0.384}
\definecolor{Fpink1}{rgb}{1.000,0.710,0.773}
\definecolor{Fpink2}{rgb}{0.933,0.663,0.722}
\definecolor{Fpink3}{rgb}{0.804,0.569,0.620}
\definecolor{Fpink4}{rgb}{0.545,0.388,0.424}
\definecolor{FLightPink1}{rgb}{1.000,0.682,0.725}
\definecolor{FLightPink2}{rgb}{0.933,0.635,0.678}
\definecolor{FLightPink3}{rgb}{0.804,0.549,0.584}
\definecolor{FLightPink4}{rgb}{0.545,0.373,0.396}
\definecolor{FPaleVioletRed1}{rgb}{1.000,0.510,0.671}
\definecolor{FPaleVioletRed2}{rgb}{0.933,0.475,0.624}
\definecolor{FPaleVioletRed3}{rgb}{0.804,0.408,0.537}
\definecolor{FPaleVioletRed4}{rgb}{0.545,0.278,0.365}
\definecolor{Fmaroon1}{rgb}{1.000,0.204,0.702}
\definecolor{Fmaroon2}{rgb}{0.933,0.188,0.655}
\definecolor{Fmaroon3}{rgb}{0.804,0.161,0.565}
\definecolor{Fmaroon4}{rgb}{0.545,0.110,0.384}
\definecolor{FVioletRed1}{rgb}{1.000,0.243,0.588}
\definecolor{FVioletRed2}{rgb}{0.933,0.227,0.549}
\definecolor{FVioletRed3}{rgb}{0.804,0.196,0.471}
\definecolor{FVioletRed4}{rgb}{0.545,0.133,0.322}
\definecolor{Fmagenta1}{rgb}{1.000,0.000,1.000}
\definecolor{Fmagenta2}{rgb}{0.933,0.000,0.933}
\definecolor{Fmagenta3}{rgb}{0.804,0.000,0.804}
\definecolor{Fmagenta4}{rgb}{0.545,0.000,0.545}
\definecolor{Forchid1}{rgb}{1.000,0.514,0.980}
\definecolor{Forchid2}{rgb}{0.933,0.478,0.914}
\definecolor{Forchid3}{rgb}{0.804,0.412,0.788}
\definecolor{Forchid4}{rgb}{0.545,0.278,0.537}
\definecolor{Fplum1}{rgb}{1.000,0.733,1.000}
\definecolor{Fplum2}{rgb}{0.933,0.682,0.933}
\definecolor{Fplum3}{rgb}{0.804,0.588,0.804}
\definecolor{Fplum4}{rgb}{0.545,0.400,0.545}
\definecolor{FMediumOrchid1}{rgb}{0.878,0.400,1.000}
\definecolor{FMediumOrchid2}{rgb}{0.820,0.373,0.933}
\definecolor{FMediumOrchid3}{rgb}{0.706,0.322,0.804}
\definecolor{FMediumOrchid4}{rgb}{0.478,0.216,0.545}
\definecolor{FDarkOrchid1}{rgb}{0.749,0.243,1.000}
\definecolor{FDarkOrchid2}{rgb}{0.698,0.227,0.933}
\definecolor{FDarkOrchid3}{rgb}{0.604,0.196,0.804}
\definecolor{FDarkOrchid4}{rgb}{0.408,0.133,0.545}
\definecolor{Fpurple1}{rgb}{0.608,0.188,1.000}
\definecolor{Fpurple2}{rgb}{0.569,0.173,0.933}
\definecolor{Fpurple3}{rgb}{0.490,0.149,0.804}
\definecolor{Fpurple4}{rgb}{0.333,0.102,0.545}
\definecolor{FMediumPurple1}{rgb}{0.671,0.510,1.000}
\definecolor{FMediumPurple2}{rgb}{0.624,0.475,0.933}
\definecolor{FMediumPurple3}{rgb}{0.537,0.408,0.804}
\definecolor{FMediumPurple4}{rgb}{0.365,0.278,0.545}
\definecolor{Fthistle1}{rgb}{1.000,0.882,1.000}
\definecolor{Fthistle2}{rgb}{0.933,0.824,0.933}
\definecolor{Fthistle3}{rgb}{0.804,0.710,0.804}
\definecolor{Fthistle4}{rgb}{0.545,0.482,0.545}
\definecolor{FDarkGrey}{rgb}{0.663,0.663,0.663}
\definecolor{FDarkGray}{rgb}{0.663,0.663,0.663}
\definecolor{FDarkBlue}{rgb}{0.000,0.000,0.545}
\definecolor{FDarkCyan}{rgb}{0.000,0.545,0.545}
\definecolor{FDarkMagenta}{rgb}{0.545,0.000,0.545}
\definecolor{FDarkRed}{rgb}{0.545,0.000,0.000}
\definecolor{FLightGreen}{rgb}{0.565,0.933,0.565}

\newcommand{\thetaprop}{\theta^{\textrm{\tiny\rm prop}}}

\RRdate{July 2011}
\RRauthor{% les auteurs
Fabien Campillo\thanks{\protect\url{Fabien.Campillo@inria.fr} --- 
           Project--Team MODEMIC, INRIA/INRA, UMR MISTEA, b\^at. 29, 2 place Viala, 34060 Montpellier cedex 06, France.}%
  \and
Dominique Hervé\thanks{\protect\url{Dominique.herve@ird.fr} --- MEM, University of Fianarantsoa \&\ IRD (UMR GRED), ENI, BP1487, 301 Fianarantsoa, Madagascar.}%
  \and
Angelo Raherinirina\thanks{University of Fianarantsoa BP 1264, Andrainjato, 301 Fianarantsoa, Madagascar}%
  \and
Rivo Rakotozafy\thanks{\protect\url{rrakotozafy@yahoo.fr} ---  University of Fianarantsoa BP 1264, Andrainjato, 301 Fianarantsoa, Madagascar}%
}
\authorhead{Campillo, Hervé, Raherinirina, Rakotozafy}
\RRnote {This works was partially supported by the LIRIMA/INRIA program and by the Cooperation and Cultural Action Office (SCAC) of the France Embassy in Madagascar.}

\RRtitle{Un modèle markovien de dynamique d'usage des terres}
\RRetitle{A Markov model of land use dynamics}
\titlehead{A Markov model of land use dynamics}
\RRabstract{The application of the Markov chain to modeling agricultural succession is well known.
In most cases, the main problem is the inference of the model, i.e. the estimation of the transition matrix. In this work we present  methods to estimate the transition matrix from historical observations. In addition to the estimator of maximum likelihood (MLE), we also consider the Bayes estimator associated with the  Jeffreys prior. This Bayes estimator  will  be approximated by a  Markov chain Monte Carlo (MCMC) method. We also propose  a method  based on the sojourn time to test the adequation  of  Markov chain model to the dataset.}

\RRkeyword{Markov model, Markov chain Monte Carlo, Jeffreys prior, land use dynamics.}

\RRresume{Les chaînes de Markov sont depuis longtemps utilisées en modélisation de la dynamique d'usage des terres. Dans la plupart des cas, se pose le problème  de l'inférence du modèle, c'est à dire de la construction de la matrice de transition qui dirige la dynamique de succession. Nous présentons dans cet article des  méthodes pour estimer cette matrice à partir d'un historique d'observations. En plus de  l'estimateur du maximum de vraisemblance (EMV), nous consi\-dérons l'estimateur bayésien associé à la loi a priori non informative de Jeffreys. Cet estimateur bayésien sera approché  par une méthode de Monte Carlo par chaîne de Markov (MCMC). Nous étudions également l'adéquation entre les temps de séjour, en un état, constatés dans les données et leur estimation par le modèle de Markov.}

\RRmotcle{Modèle de Markov, Monte Carlo par Chaîne de Markov,  loi a priori de Jeffreys, dynamique d'usage des terres.}

\RRprojet{MODEMIC}
\RRthemeProj{mere} % theme du projet Apics
\RCSophia % Sophia Antipolis M\'editerran\'ee
\begin{document}
%%%%%%%%%%%%%%%%%%%%%%%%%%%%%%%%%%%%%%%%%%%%%%%%%%%%%%%%%%%%%%%%%%%%%%%%%%%%%%%%%%%%
%%%%%%%%%%%%%%%%%%%%%%%%%%%%%%%%%%%%%%%%%%%%%%%%%%%%%%%%%%%%%%%%%%%%%%%%%%%%%%%%%%%%
%%%%%%%%%%%%%%%%%%%%%%%%%%%%%%%%%%%%%%%%%%%%%%%%%%%%%%%%%%%%%%%%%%%%%%%%%%%%%%%%%%%%
 
\RRNo{7670}
 
\makeRR   % cas d'un rapport de recherche
%% \makeRT % cas d'un rapport technique.
%% a partir d'ici, chacun fait comme il le souhaite

\clearpage
\tableofcontents
\mbox{}
\cleardoublepage
\mbox{}
\cleardoublepage

%%%%%%%%%%%%%%%%%%%%%%%%%%%%%%%%%%%%%%%%%%%%%%%%%%%%%%%%%%%%%%%%%%%%%%%%%%%%
%%%%%%%%%%%%%%%%%%%%%%%%%%%%%%%%%%%%%%%%%%%%%%%%%%%%%%%%%%%%%%%%%%%%%%%%%%%%
\section{Introduction}
%%%%%%%%%%%%%%%%%%%%%%%%%%%%%%%%%%%%%%%%%%%%%%%%%%%%%%%%%%%%%%%%%%%%%%%%%%%%
%%%%%%%%%%%%%%%%%%%%%%%%%%%%%%%%%%%%%%%%%%%%%%%%%%%%%%%%%%%%%%%%%%%%%%%%%%%%

Population pressure is one of the major causes of deforestation in tropical countries. In the region of Fianarantsoa (Madagascar), two national parks Ranomafana and Andringitra are connected by a forest corridor, which is of critical importance to maintain the regional biodiversity. The need for cultivated land pushes people to encroach on the corridor to look for swallows to be converted into paddy fields, and then to clear slope forested parcels for cultivation. Once swallows are all converted in paddy fields, the dynamic of slash and burn cultivation is clearly opposed to the dynamic of forest conservation and regeneration. To reconciling forest conservation with agricultural production, it is important to understand and model the dynamic of post-forest land use of these parcels. We will use a first data set developed by IRD, in the western edge of the corridor, consisting of the annual state of 43 parcels initially in forest, during 22 years since first clearing (Figure \ref{fig.data}). Each parcel can take four possible states:  forest ($F$), annual crop ($C$), fallow ($J$), perennial crop ($B$).

The use of Markovian approaches to model land-use transitions and vegetation successions is widespread \cite{tucker2004a,usher1979a,waggoner1970a}. The success of these approaches is explained by the fact that agro-ecological dynamics are often represented as discrete succession of a finite number of states, each one with its holding time. Both agronomists and ecologists, in dialog, actually fail in predicting the future succession of these states, knowing the previous land-use history. They ask the mathematicians for detailing the characteristics of these dynamics and defining how to pilot them. The construction, manipulation and simulation of such models are fairly easy. The transition probabilities of the Markovian model are estimated from observed data. The classical reference \cite{anderson-tw1957a} proposes the maximum likelihood method to estimate the transition probabilities of a Markov chain. An alternative is to consider Bayesian estimators \cite{minjesung2007a,meshkani1992a}. In this paper we explore and test several modeling tools, Markov chain, Bayesian estimation and MCMC procedure to better fit with the actual data. These results are needed by agro-ecologists who try to model the land use dynamics, at a parcel scale.

\medskip

The model is introduced in Section \ref{sec.model}, then the maximum likelihood estimator and the Bayesian estimator are presented in Sections \ref{sec.mle}
and  \ref{sec.bayes} respectively. These estimators are applied to simulated data in Section \ref{sec.simulation} and to the real data set in Section \ref{sec.application}. The Markov model is evaluated in Section \ref{sec.evalution}. Conclusion and perspectives are drawn in
Section \ref{sec.conlusion}.

%----------------
\begin{figure}
\begin{center}
\begin{minipage}{12.5cm}
{\scriptsize\color{darkgreen}
\begin{verbatim}
    parcel number
    1 2 3 4 ... 
      
y  0   F F F F F F F F F F F F F F F F F F F F F F F F F F F F F F F F F F F F F F F F F F F
e  1   F F F F F F F F F F F F F F F F F F F F F F F F F F F F F F F F F F F F F F F F F F C
a  2   F F F F F F F F F F F F F F F F F F F F F F F F F F F F F F F F F F F F F F F F F F C
r  3   F F F F F F F F F F F F F F F F F F F F F F F C C C C C C C C C F F F F F F F F F F C
   4   F F F F F F F F F F F F F F F F F F F F F F F C C C C C C C C J F F F F F F F F F F B
   .   F F F F F F F F F F F F F F F F F F F F F F F C C C C C C J J C F F F F F F F F F F B
   .   F F F F F F F F F F F F F F F F F F F F F F F C C J C C J J C C F F F F F F F F F F B
   .   F F F F F F F F F F F F F F F F F F F F F F F C C J C C J C C J F F F F F F F F F F B
       F F F F F F F F F F F F F F F F F F F F F F F C C J C J C J J J F F F F F F F F F F B
       F F F F F F F F F F F F F F F F F F F F F F F J J J J J C J J J F F F F F F F F F F B
       F F F F F F F F F F F F F F F F F F F F F F F J J J J C C J C C F F F F F F F F F F B
       F F F F F F F F F F F F F F F F F F F F F C C J J J J C J J C C F F F F F F F F F F B
       F F F F F F F F F F F F F F F F F F F F F C C J J J C C J C C J F F F F F F F F F F B
       F F F F F F F C F F F F F F F F F F F F F C C J J J C C C C J J F F F F F F F F F F B
       F F F F C C C C C C C F F F F F F F F F F C C J J J C C J J C J F F F F F F F F F F B
       C F F F C C C C C C C C C C C C C C C C C C C C J J C J J J C C C C C C C C C J C C B
       C C C C C C C J C C C C C C J J J C J C J J J C J J C J J J C C C C C C C C C J C C B
       J C C C C C C J C J C J C C J J J C J C C C J J C C J C C J J J C C C C C J J J C C B
       J C C C C C B C C C C J C C J J J C J J C C C J C C J C C J J J C C C C C J C J C C B
       J C C C C C B C C C C J C C J J J J J J C C C J J J C J C J C J C J C C C C C J C C B
       C B C C C J B C J J J J J J C C C J J J C J C C J C C J C C C C J C J C C J C J C C B
       C B C C C C B C J C C J J J C C J J J J J C J J C C C J C C J C C C C J C C J J C C B
\end{verbatim}
}
\end{minipage}
\end{center}
\caption{\it Annual states $(e_{n}^{p})_{n=0:21}^{p=1:43}$ corresponding to 43 parcels and 22 years. These parcels are located
on the slopes and lowlands on the edge of the forest corridor of Ranomafana-Andringitra, Madagascar. The states are: 
forest ($F$), 
fallow ($J$),
annual crop ($C$), 
perennial crop ($B$).} 
\label{fig.data}
\end{figure}
%----------------

%%%%%%%%%%%%%%%%%%%%%%%%%%%%%%%%%%%%%%%%%%%%%%%%%%%%%%%%%%%%%%%%%%%%%%%%%%%%
%%%%%%%%%%%%%%%%%%%%%%%%%%%%%%%%%%%%%%%%%%%%%%%%%%%%%%%%%%%%%%%%%%%%%%%%%%%%
\section{The model}
\label{sec.model}
%%%%%%%%%%%%%%%%%%%%%%%%%%%%%%%%%%%%%%%%%%%%%%%%%%%%%%%%%%%%%%%%%%%%%%%%%%%%
%%%%%%%%%%%%%%%%%%%%%%%%%%%%%%%%%%%%%%%%%%%%%%%%%%%%%%%%%%%%%%%%%%%%%%%%%%%%

We make the following hypothesis:
\begin{itemize}
\item[($H_{1}$)] \textit{The dynamics of the parcels are independent and identical.}
\end{itemize}
This means that  $(e_{n}^{p})_{n=0:21}^{p=1:43}$ are  43 independent realizations of a same process $(X_n)_{n=0:21}$. This assumption is not realistic as the dynamics of a given parcel depends on:
\begin{itemize}
\item farmer decisions;
\item exposition, slope and distance from the forest, that means properties of the same plot;
\item neighboring parcels.
\end{itemize}
This assumption, however will lead to a simple model.

We also suppose that:
\begin{itemize}
\item[($H_{2}$)] \textit{The process $(X_n)_{n=0\cdots 21}$ is Markovian and time-homogeneous.}
\end{itemize}
The homogeneity assumption is also simplistic but we assume that the transition law of parcels will poorly varied during this 22 year period.

Finally we suppose that:
\begin{itemize}
 \item[($H_{3}$)] \textit{The initial state is $F$.}
\end{itemize}

%--------------------------------
\begin{figure}
 \begin{center}
\includegraphics[width=8cm]{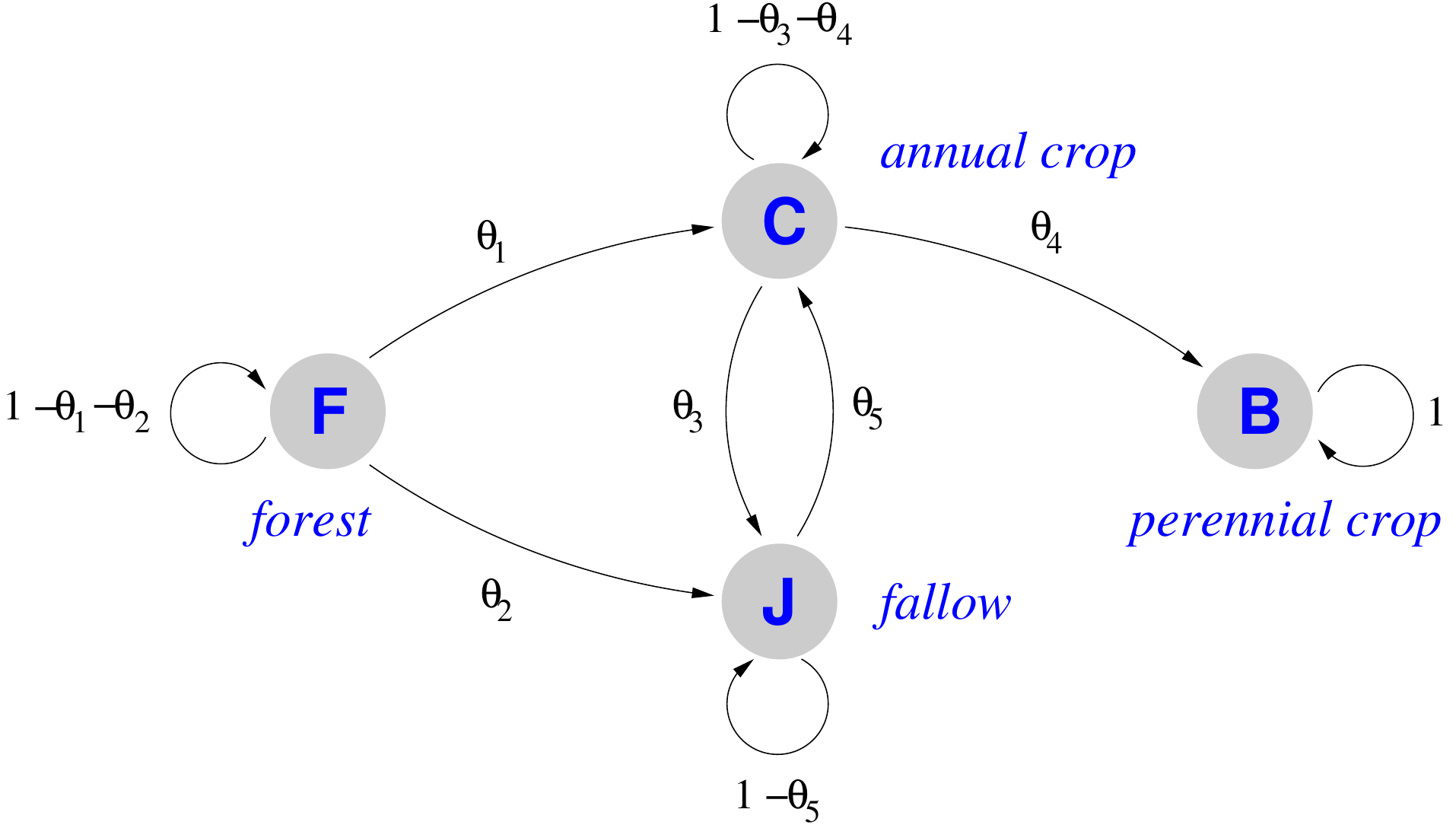}
\caption{\it Four states Markov chain diagram: forest ($F$), 
annual crop ($C$), 
fallow ($J$),
perennial crop ($B$); $F$ is the initial state and  $B$ is an absorbing state.}
\end{center}
\label{fig2}
\end{figure}
%--------------------------------

These hypotheses lead to a model $X=(X_n^p)^{p=1:P}_{n=0:N-1}$, $P=43$, 
$N=22$ where $(X_n^p)_{n=0\cdots N-1}$ are $P$ independent Markov chains, with initial law $\delta_F$  and transition matrix $Q$ of size $4 \times 4$. The state space is:
\[
  E \eqdef \{ F, C, J, B\}\,.
\]
Hence:
\begin{subequations}
\label{eq.model}
\begin{align}
\label {eq1}
 \nonumber
 \P(X_{0\:N-1}^{1:P}= e_{0:N-1}^{1:P}) 
  &= 
  \prod_{p=1}^{P} 
     \P(X_{0:N-1}^p=e_{0:N-1}^p) 
\\
  &=
  \prod_{p=1}^{P}\,\delta_F\,(e_0^p)\;Q(e_0^p, e_1^p)
   \cdots 
 Q(e_{N-2}^p,e_{N-1}^p)
\end{align}
for all $ e_n^p\in E$.

Some transitions are not observed at least in the time scale considered here: once the parcel leaves the state ``forest'' by first clearing, it cannot come back; and similarly when it reaches the state ``perennial crop'', it stays there during the sample time, that means permanently in this model.
\begin{itemize}
 \item[($H_{4}$)] \textit{The transitions $C\to F$, $J\to F$, $B\to F$, $B\to C$, $B\to J$, $J\to B $ do not exist in the model, all other transitions are possible.}
\end{itemize}
In particular: once the parcel leaves the state ``forest'', it cannot come back; when it reaches the state  ``perennial crop'', it stays there permanently.

This hypothesis implies that \fenumi\ the realistic transitions, which are not observed during the considered time scale, $J\to F$, $B\to C$, do not exist in the model, \fenumii\ the unrealistic transitions $C\to F$, $B\to F$, $B\to J$ and $J\to B$ do not exist, \fenumiii\   the state $F$ is transient (more precisely when the chain leaves the state $F$ it will never come back to that state), \fenumiv\  the state $B$ is absorbing.

\bigskip

To summarize we consider a transition matrix $Q$ of the form:
\begin{equation}
\label{eq2}
Q = \left(
         \begin{array}{cccc}
             1-\theta_1-\theta_2 & \theta_1 & \theta_2 & 0 \\
             0 & 1-\theta_3-\theta_4 & \theta_3 & \theta_4 \\
            0 & \theta_5 & 1-\theta_5 & 0\\
             0 & 0 & 0 & 1 \\
         \end{array}
      \right ),
\end{equation}
\end{subequations}
that corresponds to Figure \ref{fig2}, and depends on a 5-dimensional parameter:
\[
 \theta=(\theta_1, \theta_2, \theta_3, \theta_4, \theta_5)
\]
belonging to the set:
\begin{align}
\label{eq.Theta}
   \Theta 
   \eqdef 
   \big\{ 
     \theta \, \in [0,\,1]^5
     \,;\,\;
     \theta_1 +  \theta_2\leq 1 ,\; 
     \theta_3 +  \theta_4\leq 1 
    \big\}\,.
\end{align}
Let $\P_{\theta}$ denotes the probability under which the Markov chain admit $Q$ with parameter $\theta$ as a transition matrix.

%%%%%%%%%%%%%%%%%%%%%%%%%%%%%%%%%%%%%%%%%%%%%%%%%%%%%%%%%%%%%%%%%%%%%%%%%%
%%%%%%%%%%%%%%%%%%%%%%%%%%%%%%%%%%%%%%%%%%%%%%%%%%%%%%%%%%%%%%%%%%%%%%%%%%
\section{Maximum likelihood estimation}
\label{sec.mle}
%%%%%%%%%%%%%%%%%%%%%%%%%%%%%%%%%%%%%%%%%%%%%%%%%%%%%%%%%%%%%%%%%%%%%%%%%%
%%%%%%%%%%%%%%%%%%%%%%%%%%%%%%%%%%%%%%%%%%%%%%%%%%%%%%%%%%%%%%%%%%%%%%%%%%

We recall the classical results of Anderson-Goodman \cite{anderson-tw1957a} to compute the MLE of the matrix $Q$.

The likelihood function associated with $\{\P_{\theta};\theta\in\Theta\}$ is:
\begin{align*}
   L(\theta) 
   \eqdef \P_{\theta}(X_{0:N-1}^{1:P} = e_{0:N-1}^{1:P})
   = \prod_{p=1}^{P}\,\delta_F\,(e_0^p)\;Q(e_0^p, e_1^p)\cdots Q(e_{N-2}^p,e_{N-1}^p)
\end{align*}
where $Q$ is defined by \eqref{eq2}, for any $e_{0:N-1}^{1:P}\in E^{P\times N}$. Let $n_{ee'}^p$ be the number of transitions from state $e$ to state $e'$ for a parcel $p$ in $X$:
\begin{align}
\label{eq3}
  n_{ee'}^p 
  \eqdef 
  \sum_{n=1}^{N-1} \indic_{\{X_{n-1}^p=e\}}\; \indic_{\{X_n^p = e'\}}
  \qquad \forall e,e' \in E
\end{align}
and 
$n_{ee'}$ be the total of number of transitions from state $e$ to state $e'$ ($e , e' \in E$):
\begin{align} 
\label{eq4}
  n_{ee'} \eqdef \sum_{p=1}^{P} n_{ee'}^p\,.
\end{align}
According to \eqref{eq2}:
\begin{align*}
  &
  L(\theta)
  =\prod_{p=1}^P Q(1,1)^{n_{FF}^p}\; Q(1,2)^{n_{FC}^p} \, Q(1,3)^{n_{FJ}^p}
  \, Q(2,2)^{n_{CC}^p}\,Q(2,3)^{n_{CJ}^p}
\\[-0.7em]
  &
  \hspace{4cm} Q(2,4)^{n_{CB}^p}\, Q(3,2)^{n_{JC}^p}
  Q(3,3)^{n_{JJ}^p}
\\
  &
  \quad
  = \prod_{p=1}^{P} (1-\theta_1 - \theta_2)^{n_{FF}^p}\, 
            \theta_1^{n_{FC}^p}\, \theta_2^{n_{FJ}^p}\,(1- \theta_3 - \theta_4)^{n_{CC}^p}
    \,\theta_3^{n_{CJ}^p}
    \,\theta_4^{n_{CB}^p} \, (1 - \theta_5 )^{n_{JJ}^p}\,\theta_5^{n_{JC}^p}
\end{align*}
and from \eqref{eq4}:
\begin{align*}
  &
  L(\theta)
  \eqdef (1-\theta_1 - \theta_2)^{n_{FF}}\, \theta_1^{n_{FC}}\, \theta_2^{n_{FJ}}\,
    (1- \theta_3 - \theta_4)^{n_{CC}}\, \theta_3^{n_{CJ}}\,
  \theta_4^{n_{CB}} \, (1 - \theta_5 )^{n_{JJ}}\,
     \theta_5^{n_{JC}}
\end{align*}
so that the log-likelihood function reads:
\begin{align}
\nonumber
  l(\theta)
  &=  
  n_{FF}  \log(1-\theta_1 -\theta_2) + n_{FC}\log(\theta_1) + n_{FJ}\log(\theta_2) 
\\
\nonumber
  & 
  \qquad \qquad 
  + n_{CC}\log(1-\theta_3 - \theta_4) + n_{CJ}\log(\theta_3) + n_{CB}\log(\theta_4)
\\
\label{eq6}
  & 
  \qquad \qquad
  + n_{JJ} \log(1-\theta_5) + n_{JC}\log(\theta_5) \,.
\end{align}
The MLE $\hat\theta$ is solution of $\partial l(\theta)/\partial \theta|_{\theta=\hat\theta}=0 $, that is:
\begin{align*}
  \frac{\partial l(\theta)}{\partial \theta_1}
  &=
  -\frac{n_{FF}}{(1-\theta_1 - \theta_2)} +  \frac{n_{FC}}{\theta_1} = 0
  \,,
&
  \frac{\partial l(\theta)}{\partial \theta_2} 
  &= 
  -\frac{n_{FF}}{(1-\theta_1 - \theta_2)} + \frac{n_{FJ}}{\theta_2} = 0
  \,,
\\
  \frac{\partial l(\theta)}{\partial \theta_3} 
  &= 
  -\frac{n_{CC}}{(1-\theta_3 - \theta_4)} + \frac{n_{CJ}}{\theta_3} = 0
  \,,
&
  \frac{\partial l(\theta)}{\partial \theta_4} 
  &= 
  -\frac{n_{CC}}{(1-\theta_3 - \theta_4)} + \frac{n_{CB}}{\theta_4} = 0
  \,, 
\\
  \frac{\partial l(\theta)}{\partial \theta_5} 
  &= 
  -\frac{n_{JJ}}{(1-\theta_5 )} + \frac{n_{JC}}{\theta_5} = 0
  \,. 
\end{align*}
We get:
\begin{align*}
  \hat{\theta}_1 
  &=\frac{n_{FC}}{n_{FF} + n_{FC} + n_{FJ}}\,,
  &
  \hat{\theta}_2
  &=\frac{n_{FJ}}{n_{FF} + n_{FC} + n_{FJ}}\,,
\\
  \hat{\theta}_3 
  &=\frac{n_{CJ}}{n_{CC} + n_{CF} + n_{CJ}}\,,
  &
  \hat{\theta}_4
  &=\frac{n_{CB}}{n_{CC} + n_{CJ} + n_{CB}}\,,
\\
  \hat{\theta}_5
  &=\frac{n_{JC}}{n_{JC} + n_{JJ} }\,.
\end{align*}

%%%%%%%%%%%%%%%%%%%%%%%%%%%%%%%%%%%%%%%%%%%%%%%%%%%%%%%%%%%%%%%%%%%%%%%%%%
%%%%%%%%%%%%%%%%%%%%%%%%%%%%%%%%%%%%%%%%%%%%%%%%%%%%%%%%%%%%%%%%%%%%%%%%%%
\section{Bayesian estimation}
\label{sec.bayes}
%%%%%%%%%%%%%%%%%%%%%%%%%%%%%%%%%%%%%%%%%%%%%%%%%%%%%%%%%%%%%%%%%%%%%%%%%%
%%%%%%%%%%%%%%%%%%%%%%%%%%%%%%%%%%%%%%%%%%%%%%%%%%%%%%%%%%%%%%%%%%%%%%%%%%

We suppose that an a priori distribution law $\pi (\theta)$ on the parameter $\theta$ is given. According to the Bayes rule, the a posteriori distribution law $\pi (\theta)$ on $\theta$ given the observations $X$ is:
\begin{align}
\label{eq.posterior}
  \tilde{\pi}(\theta) \propto L(\theta) \; \pi (\theta)
\end{align}
where $L(\theta)$ is the likelihood function. The Bayes estimator $\tilde \theta $ of the parameter $\theta$ is the mean of the a posteriori distribution:
\begin{align}
\label{eq.bayes.estimator}
 \tilde \theta \eqdef \int_{\Theta} \theta\,\tilde{\pi}(\theta)\,\rmd \theta\,.
\end{align}

%%%%%%%%%%%%%%%%%%%%%%%%%%%%%%%%%%%%%%%%%%%%%%%%%%%%%%%%%%%%%%%%%%%%%%%%%%
\subsection{Jeffreys prior}
%%%%%%%%%%%%%%%%%%%%%%%%%%%%%%%%%%%%%%%%%%%%%%%%%%%%%%%%%%%%%%%%%%%%%%%%%%

Numerical tests that will be performed in Section \ref{sec.simu.2.states} suggest that the Jeffreys prior is well adapted to the present situation and we introduce it now. This prior distribution (non-informative) is defined by \cite{marin2007a}:
\begin{equation}
\label{def.jeffreys}
  \pi(\theta) \propto \sqrt{\det[\II(\theta)]}
\end{equation}
where $\II(\theta)$ is the Fisher information matrix given by:
\begin{equation*}
  \II(\theta) 
  \eqdef
  \left[\E_{\theta} 
    \Big(
       -\frac{ \partial^2 l(\theta)}{ \partial \theta_k\, \partial \theta_l}
    \Big)
  \right]_{1\leq k,l\leq 5}
\end{equation*}
and $l(\theta)$ is the log-likelihood function. Hence:
\begin{equation*}
  \II(\theta) =  \left(
  \begin{array}{ccc}
      A_{1,2} &  0   &   0 \\
      0 &  A_{3,4}   &   0 \\
      0 &  0        &   a_{5} 
  \end{array}  
  \right)
\end{equation*}
with
\[
  A_{k,\ell} \eqdef -\E_{\theta} \left(
  \begin{array}{cc}
   \frac{ \partial^2 l(\theta)}{ \partial^2 \theta_k} & 
   \frac{ \partial^2 l(\theta)}{ \partial \theta_k \partial \theta_{\ell}} \\
   \frac{ \partial^2 l(\theta)}{ \partial \theta_k \partial \theta_{\ell}} &
   \frac{ \partial^2 l(\theta)}{ \partial^2 \theta_{\ell}} 
  \end{array}  
  \right)
  \,,\quad
  a_{5}\eqdef   -\E_{\theta} \left(
    \frac{ \partial^2 l(\theta)}{ \partial^2 \theta_5}
  \right)
    \,.
\]
So  $\det{\II(\theta)}= \det{A_{1,2}}\,\times\,\det{A_{3,4}}\,\times\,a_{5}$ and
\begin{subequations}
\label{eq.jeffreys}
\begin{equation}
\label{eq.jeffreys.1}
  \pi(\theta) \propto \sqrt{\det{A_{1,2}}\,\times\,\det{A_{3,4}}\,\times\,a_{5}}\,.
\end{equation}
According to \eqref{eq6}:
\begin{align*}
  \textstyle
   \frac{\partial^2 l(\theta)}{ \partial^2 \theta_1} 
   & = 
  \textstyle
            - \left( \frac{n_{FF}}{(1-\theta_1-\theta_2)^2}
              + \frac{n_{FC}}{\theta_1^2}\right)\,,
   &
  \textstyle
         \frac{\partial^2 l(\theta)}{\partial \theta_1 \,\partial
            \theta_2}
   &
  \textstyle
     = - \frac{n_{FF}}{(1-\theta_1-\theta_2)^2}\,,
\\
  \textstyle
   \frac{\partial^2 l(\theta)}{ \partial^2 \theta_2} 
   &
  \textstyle
    = 
            -\left( \frac{n_{FF}}{(1-\theta_1-\theta_2)^2} 
             + \frac{n_{FJ}}{\theta_2^2}\right)\,,
   &
  \textstyle
         \frac{\partial^2 l(\theta)}{\partial \theta_2 \,\partial
            \theta_1} 
   &
  \textstyle
   = 
   - \frac{n_{FF}}{(1-\theta_1-\theta_2)^2}\,,
\\
  \textstyle
   \frac{\partial^2 l(\theta)}{ \partial^2 \theta_3} 
   &
   = 
  \textstyle
            -\left(\frac{n_{CC}}{(1-\theta_3-\theta_4)^2} 
             + \frac{n_{CJ}}{\theta_3^2}\right)\,,
   & 
  \textstyle
         \frac{\partial^2 l(\theta)}{\partial \theta_3 \,\partial
            \theta_4}
   &=
  \textstyle
    - \frac{n_{CC}}{(1-\theta_3-\theta_4)^2}\,,
\\
  \textstyle
   \frac{\partial^2 l(\theta)}{ \partial^2 \theta_4} 
   &= 
  \textstyle
            -\left(\frac{n_{CC}}{(1-\theta_3-\theta_4)^2} 
             + \frac{n_{CB}}{\theta_4^2}\right)\,,
   &
  \textstyle
         \frac{\partial^2 l(\theta)}{\partial \theta_4 \,\partial
            \theta_3}
   &= 
  \textstyle
   - \frac{n_{CC}}{(1-\theta_3-\theta_4)^2}\,,
\\
  \textstyle
   \frac{\partial^2 l(\theta)}{ \partial^2 \theta_5}
   & = 
  \textstyle
            -\left(\frac{n_{JJ}}{(1-\theta_5)^2} 
             + \frac{n_{JC}}{\theta_5^2}\right)
\end{align*}
and
\begin{align}
\label{eq.jeffreys2a}
  \det{A_{1,2}} &=
  \textstyle
     \left( \frac{\E_\theta[n_{FF}]}{(1-\theta_1-\theta_2)^2}
        + \frac{\E_\theta[n_{FC}]}{\theta_1^2}
     \right)
     \;
     \left( \frac{\E_\theta[n_{FF}]}{(1-\theta_1-\theta_2)^2} 
        + \frac{\E_\theta[n_{FJ}]}{\theta_2^2} 
     \right) 
     -
     \left( 
        \frac{\E_\theta[n_{FF}]}{(1-\theta_1-\theta_2)^2} 
     \right)^2\,, 
\\
\label{eq.jeffreys2b}
   \det{A_{3,4}} 
    &= 
  \textstyle
     \left( \frac{\E_\theta[n_{CC}]}{(1-\theta_3-\theta_4)^2} 
        + \frac{\E_\theta[n_{CJ}]}{\theta_3^2} 
     \right)
     \;
     \left(\frac{\E_\theta[n_{CC}]}{(1-\theta_3-\theta_4)^2} 
             + \frac{\E_\theta[n_{CB}]}{\theta_4^2} 
     \right) 
     -
     \left(
        \frac{\E_\theta[n_{CC}]}{(1-\theta_3-\theta_4)^2}
     \right)^2\,,
\\
\label{eq.jeffreys2c}
 a_{5}
 &=      
  \textstyle
 \frac{\E_\theta[n_{JJ}]}{(1-\theta_5)^2} 
             + \frac{\E_\theta[n_{JC}]}{\theta_5^2} 
   \,.
\end{align}
From \eqref{eq3} and \eqref{eq4}:
\begin{align}
\nonumber
   \E_{\theta}[n_{ee'}] 
   & =  
   \sum_{p=1}^P \sum_{n=2}^N \P_{\theta}(X_n^p = e',X_{n-1}^p = e)
 \\
\nonumber
   & =  
   P\,\sum_{n=2}^N \P_{\theta}(X_n = e'|X_{n-1} = e)\,
                       \P_{\theta}(X_{n-1} = e)
\\
\nonumber
  &=  
  P\,Q(e,e') \sum_{n=2}^N \P_{\theta}(X_{n-1} = e)  
  =  
  P\,Q(e,e') \sum_{n=2}^N \left (\delta_F Q^{(n-1)} \right )_e 
\\
\label{eq.jeffreys3}
 & =  
 P\,Q(e,e') \sum_{n=2}^N [Q^{(n-1)}](F,e)
\end{align}
\end{subequations}
for all $e,e'\in E$. 

Note that \cite{assoudou2004a} proposed a more complex method to compute the Jeffreys prior distribution.

%%%%%%%%%%%%%%%%%%%%%%%%%%%%%%%%%%%%%%%%%%%%%%%%%%%%%%%%%%%%%%%%%%%%%%%%%%
\subsection{MCMC method}
%%%%%%%%%%%%%%%%%%%%%%%%%%%%%%%%%%%%%%%%%%%%%%%%%%%%%%%%%%%%%%%%%%%%%%%%%%

Although the Jeffrey prior distribution is explicit, we cannot compute analytically the corresponding Bayes estimator. We propose to use a Monte Carlo Markov chain (MCMC) method, namely a Metropolis-Hastings algorithm with a Gaussian proposal dsitribution, see  Algorithm~\ref{alg1}.

%--------------------
\begin{algorithm}
\caption{\it MCMC method: Metropolis-Hastings algorithm with a Gaussian proposal distribution. The target distribution is  $\tilde{\pi}(\theta)$ defined by \eqref{eq.posterior}, the Gaussian proposal PDF (probability density function) is $g(\cdot-\theta)$ ($g$ PDF of the $\NN(0,\sigma^2)$ distribution) where $\theta$ is the current value of the parameter.}
\label{alg1}
\begin{center}
\begin{minipage}{11cm}
\hrulefill\\[-1em]
%\mbox{}
\begin{algorithmic}
\STATE choose $\theta$
\STATE $\bar{\theta} \ot \theta$
\FOR {$k=2,3,4,\dots$}
\STATE $\epsilon\sim \NN(0,\sigma^2)$
  \STATE $\thetaprop \ot \theta + \epsilon $
  \STATE $u\sim U[0,1]$ 
  \STATE $\alpha \ot  
    \min\Big\{
      1
      , \displaystyle
      \frac{
         {\tilde{\pi}(\thetaprop)\, g(\theta-\thetaprop)}
      }{
        {\tilde{\pi}(\theta)}\,g(\thetaprop-\theta)
      }\Big\}$
\% $g $ PDF of $\NN(0,\sigma^2)$
  \IF {$u\leq \alpha$}
      \STATE $\theta \ot \thetaprop$
      \% acceptation
  %\ELSE
   %   \STATE $\theta_{k} \ot \theta_{k-1}$
 %\COMMENT  rejet
  \ENDIF
  \STATE $\displaystyle\bar{\theta} \ot \frac{k-1}{k}\,\bar{\theta} + \frac{1}{k}\,\theta$
\ENDFOR
\end{algorithmic}
\hrulefill
\end{minipage}
\end{center}
\end{algorithm}
%--------------------

%%%%%%%%%%%%%%%%%%%%%%%%%%%%%%%%%%%%%%%%%%%%%%%%%%%%%%%%%%%%%%%%%%%%%%%%%%
%%%%%%%%%%%%%%%%%%%%%%%%%%%%%%%%%%%%%%%%%%%%%%%%%%%%%%%%%%%%%%%%%%%%%%%%%%
\section{Simulation tests}
\label{sec.simulation}
%%%%%%%%%%%%%%%%%%%%%%%%%%%%%%%%%%%%%%%%%%%%%%%%%%%%%%%%%%%%%%%%%%%%%%%%%%
%%%%%%%%%%%%%%%%%%%%%%%%%%%%%%%%%%%%%%%%%%%%%%%%%%%%%%%%%%%%%%%%%%%%%%%%%%

%%%%%%%%%%%%%%%%%%%%%%%%%%%%%%%%%%%%%%%%%%%%%%%%%%%%%%%%%%%%%%%%%%%%%%%%%%
\subsection{Two states case}
\label{sec.simu.2.states}
%%%%%%%%%%%%%%%%%%%%%%%%%%%%%%%%%%%%%%%%%%%%%%%%%%%%%%%%%%%%%%%%%%%%%%%%%%

We first consider the simpler two states case $E=\{0,1\}$. It has no connection with the Markov model considered in the present work but it allows to easily compare the following different prior distributions:
\begin{enumerate}
\item the uniform distribution;
\item the beta distribution of parameter $(\frac{1}{2},\frac{1}{2})$;
\item the non-informative Jeffreys distribution.
\end{enumerate}
The Bayesian estimator is explicit for the two first priors, see Appendix A.

\begin{figure}
\centering
\includegraphics[width=4cm]{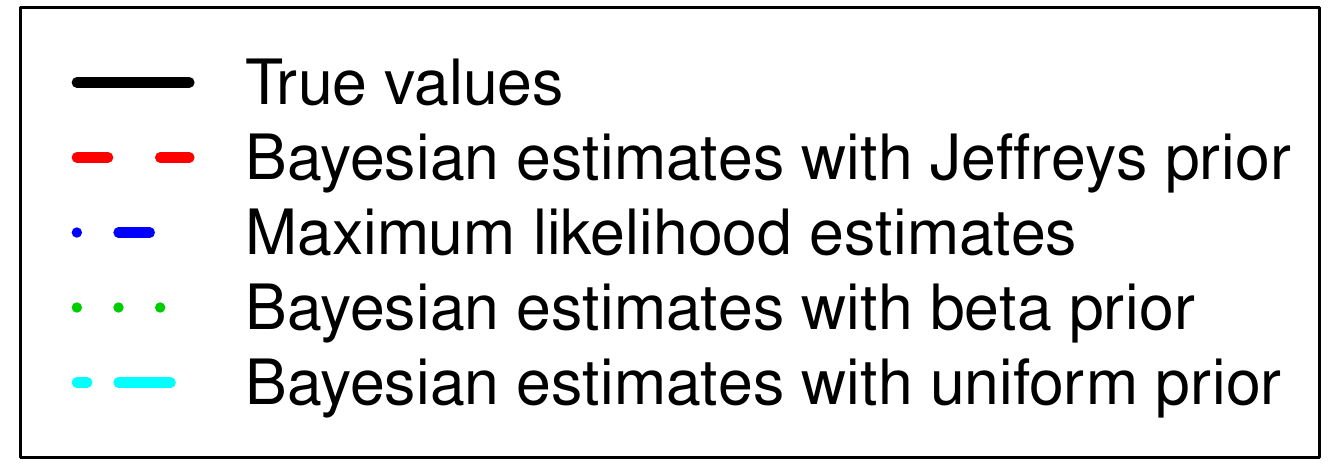}\\[0.8em]
\begin{tabular}{cc}
\includegraphics[width=8cm]{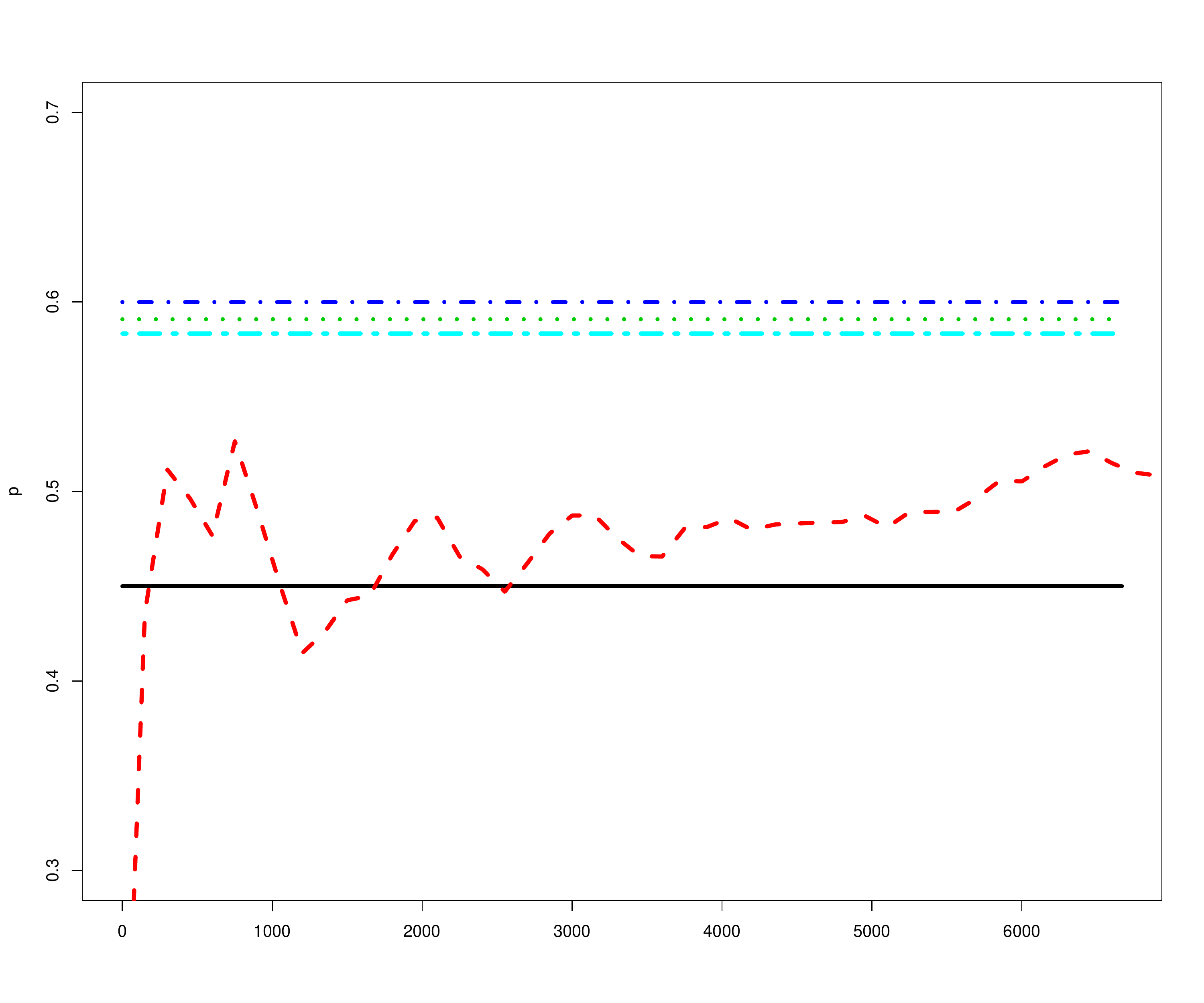}
&
\includegraphics[width=8cm]{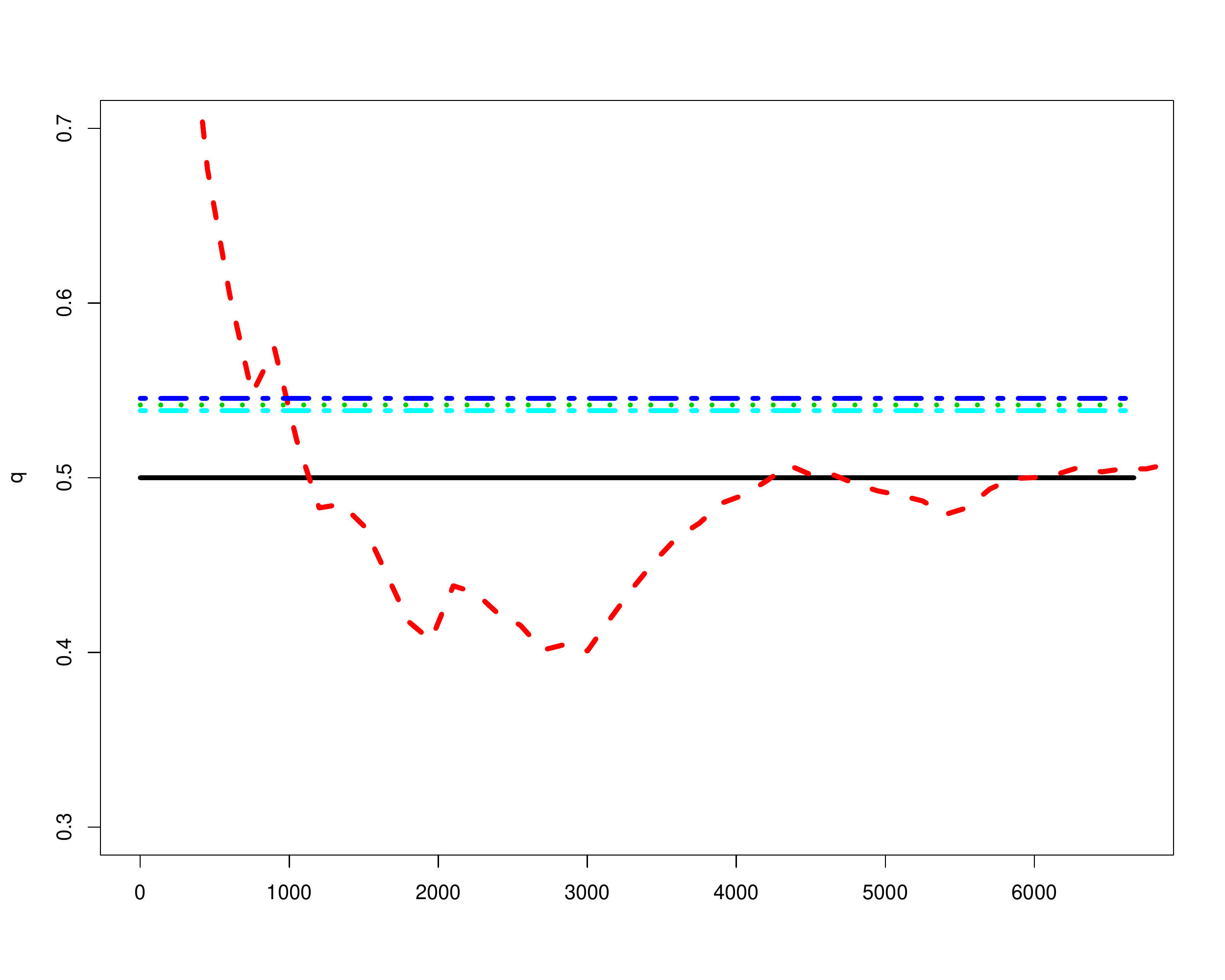}
\\
\small Estimation of $p$
&
\small Estimation of $q$
\end{tabular}
\caption{\it Two states case with paramaters $p=\P(X_{n+1}=0|X_{n}=0)$ and $q=\P(X_{n+1}=1|X_{n}=1)$: we compare the following priors: uniform, beta($1/2,1/2$) and  Jeffreys. Jeffreys prior gives better results than the two other priors.}
\label{two-state}
\end{figure}

We compare the MLE and the Bayesian estimator with uniform and beta priors, that can be explicitly computed, see Appendix A, with the Bayesian estimator with Jeffreys prior that is computed by an MCMC method that will be explained later. Results proposed in Figure \ref{two-state} tend to demonstrate that the Jeffreys prior gives better results than the two other priors.

%------------------
%------------------

%%%%%%%%%%%%%%%%%%%%%%%%%%%%%%%%%%%%%%%%%%%%%%%%%%%%%%%%%%%%%%%%%%%%%%%%%%
\subsection{Four states case}
\label{sec.simu.4.states}
%%%%%%%%%%%%%%%%%%%%%%%%%%%%%%%%%%%%%%%%%%%%%%%%%%%%%%%%%%%%%%%%%%%%%%%%%%

Before processing the real data set of Figure \ref{fig.data} with a four states Markov model, we consider a simulated case test. We aim to compare the MLE and the Bayes estimator with the Jeffreys prior.

We compute the distance between the real transition matrix $Q$ and its estimation, with the MLE or the Bayes estimator, given by the Frobenius norm:
\begin{align}
\label{eq.2-norm}
   \|A\|^2_F \eqdef \textrm{trace}(A^*A)
\end{align}
and the 2-norm:
\begin{align}
\label{eq.frobenius-norm}
  \|A\|_2 \eqdef \sqrt{\lambda_{\textrm{\tiny max}}(A^*\,A)} 
\end{align}
where $\lambda_{\textrm{\tiny max}}(A^*\,A)$ is the largest eigenvalue of the matrix $A^*\,A$.

We sample 1000 independent values $\theta^{(\ell)}$ of the parameter according to a uniform distribution on $\Theta$ defined by \eqref{eq.Theta}, that is a uniform distribution on $[0,1]^5$ with specific the constraints. For each $\ell$, we simulated data $(e_{n}^{p})_{n=0:21}^{p=1:43}$ according to the model \eqref{eq.model} that is with the transition matrix $Q_{\theta^\ell}$ defined by \eqref{eq1}. Then we compute the MLE $\hat\theta^{(\ell)}$ and the Bayes estimate $\tilde\theta^{(\ell)}$ with the Jeffreys prior. Then we compute the errors:
\begin{subequations}
\label{eq.errors}
\begin{align}
\label{eq.errors.1}
  \hat \epsilon^{\ell}
  &=\|Q_{\theta^{(\ell)}}-  Q_{\hat\theta^{(\ell)}}\|\,,
  \\
\label{eq.errors.2}
  \tilde \epsilon^{\ell}
  &=\|Q_{\theta^{(\ell)}}-  Q_{\tilde\theta^{(\ell)}}\|
\end{align}
\end{subequations}
for the two different norms. 

\begin{figure}
\centering
\includegraphics[width=5cm]{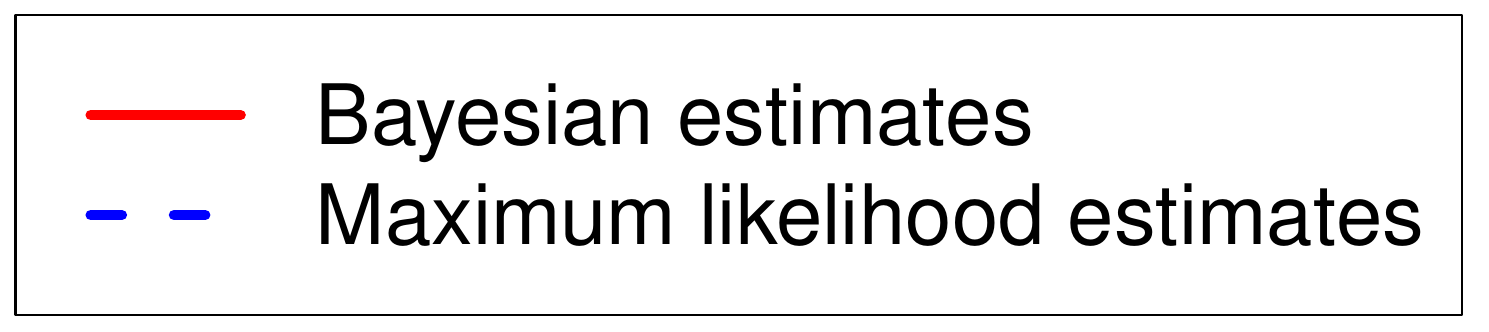}\\[0.8em]
\includegraphics[width=14cm]{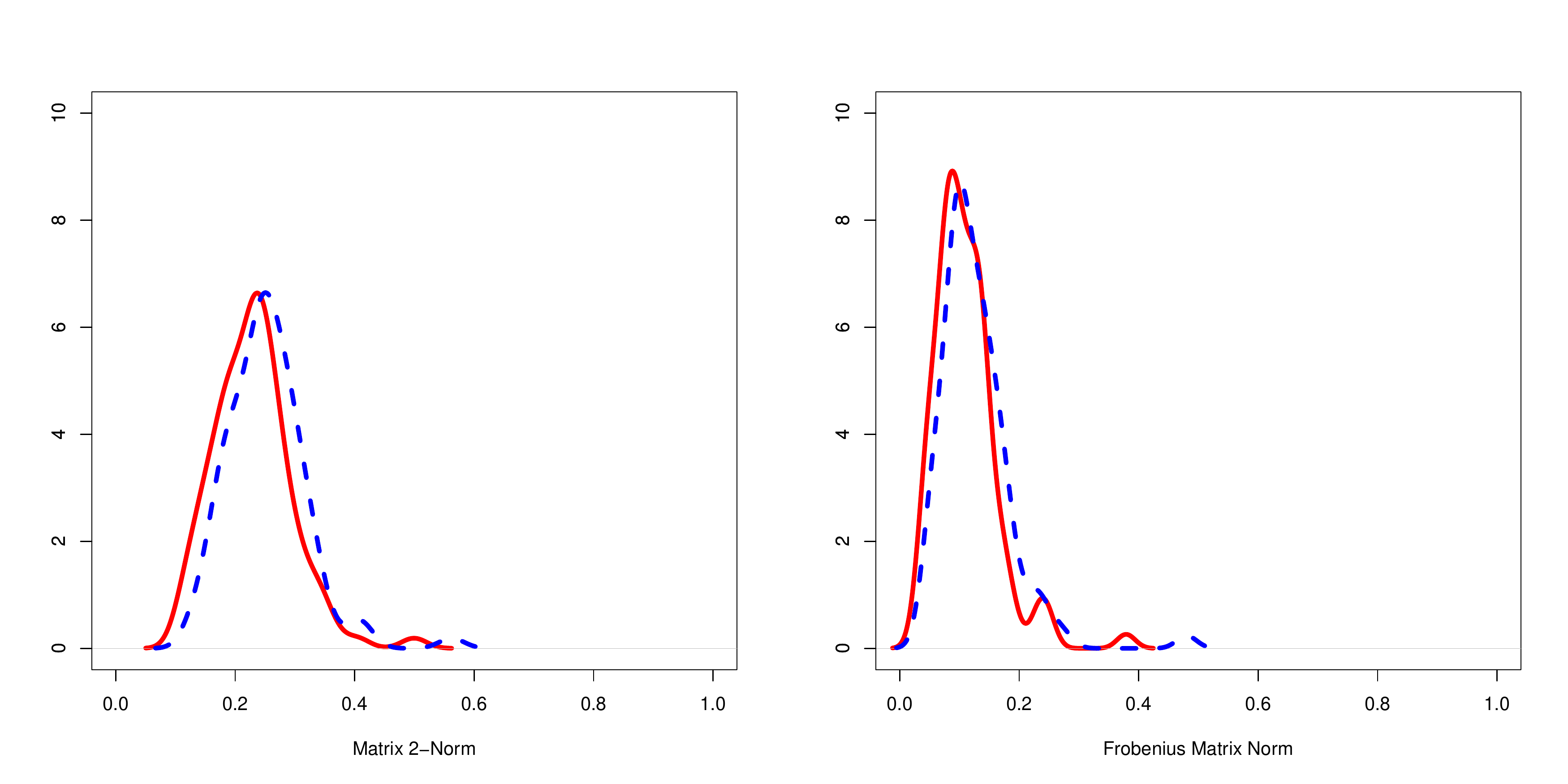}\\
\caption{\it Empirical PDF for the error terms \eqref{eq.errors} associated with the 2-norm \eqref{eq.2-norm} and the Frobenius norm \eqref{eq.frobenius-norm} based on 1000 simulation of the parameter $\theta$.}
\label{fig.sim.4state}
\end{figure}

In Figure \ref{fig.sim.4state}
 we plotted the empirical distribution of the errors $\hat \epsilon^{\ell}$ and $\tilde \epsilon^{\ell}$,  $\ell=1\cdots 1000$, for the two different norms. We see that the Bayes estimator give slightly better results than the MLE.

%------------------
%------------------

%%%%%%%%%%%%%%%%%%%%%%%%%%%%%%%%%%%%%%%%%%%%%%%%%%%%%%%%%%%%%%%%%%%%%%%%%%
%%%%%%%%%%%%%%%%%%%%%%%%%%%%%%%%%%%%%%%%%%%%%%%%%%%%%%%%%%%%%%%%%%%%%%%%%%
\section{Application to the real data set}
\label{sec.application}
%%%%%%%%%%%%%%%%%%%%%%%%%%%%%%%%%%%%%%%%%%%%%%%%%%%%%%%%%%%%%%%%%%%%%%%%%%
%%%%%%%%%%%%%%%%%%%%%%%%%%%%%%%%%%%%%%%%%%%%%%%%%%%%%%%%%%%%%%%%%%%%%%%%%%

\begin{figure}
\centering
\begin{tabular}{c@{\hskip4em}c}
\multicolumn{2}{c}
{\includegraphics[width=4.5cm]{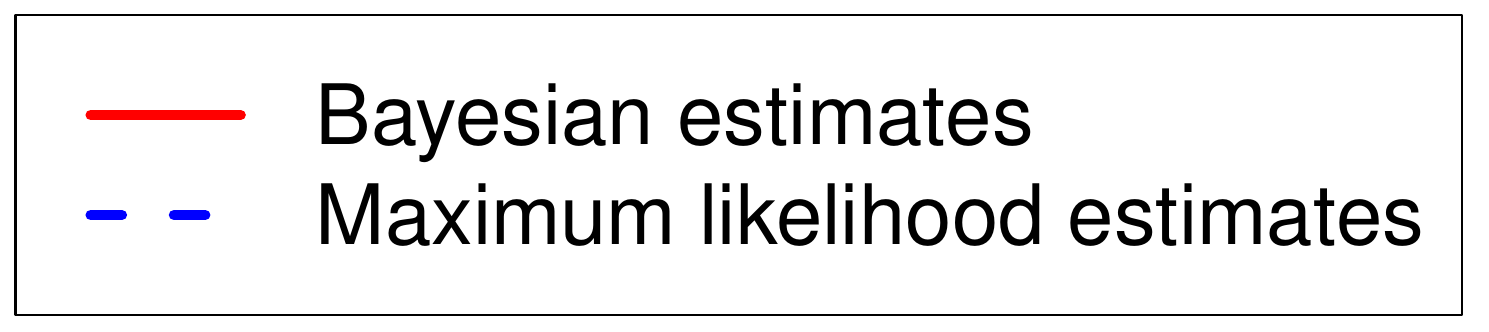}}
\\[0.5em]
\rotatebox{90}{\hskip 3.5em \scriptsize MCMC iterations}
\includegraphics[width=5cm]{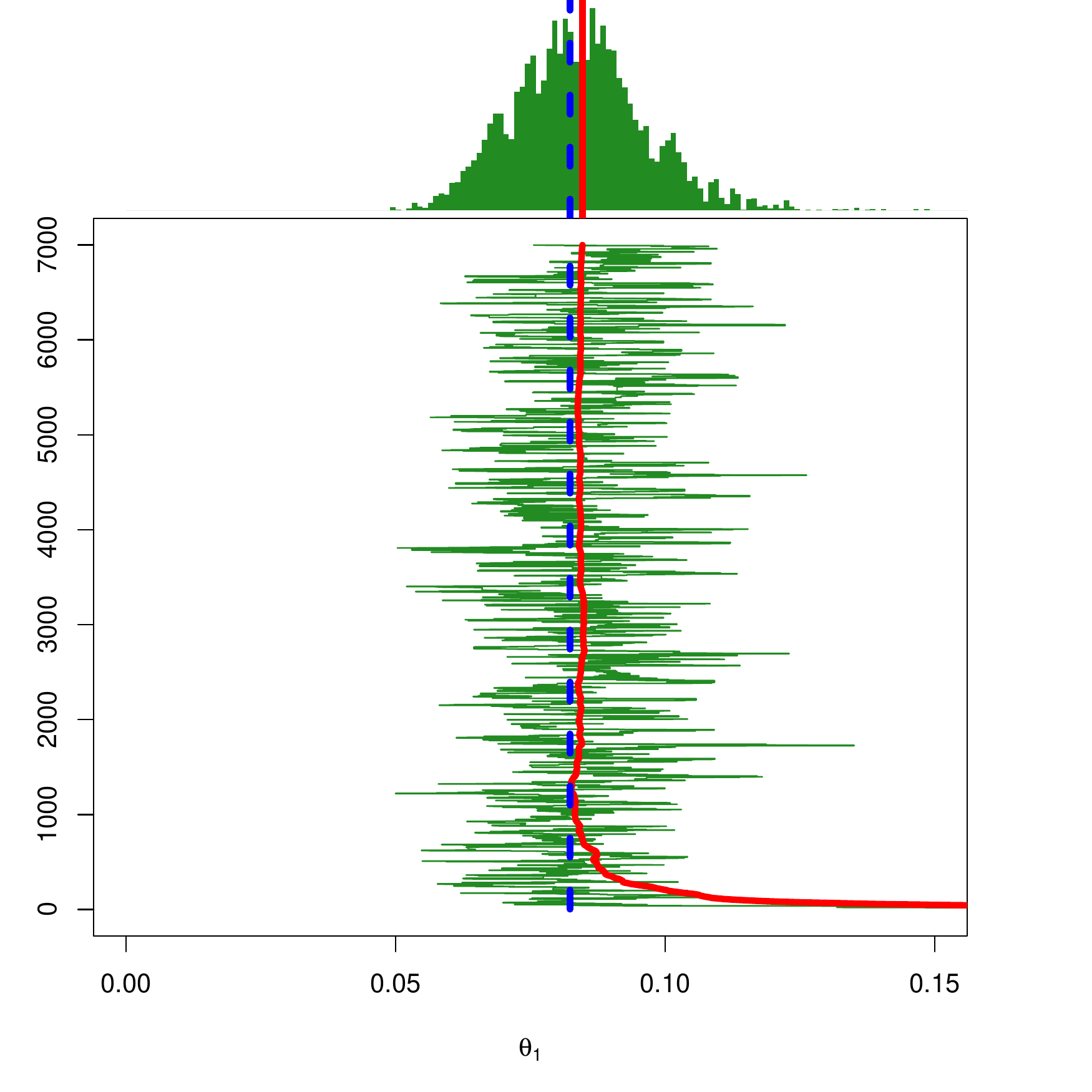}
&
\rotatebox{90}{\hskip 3.5em \scriptsize MCMC iterations}
\includegraphics[width=5cm]{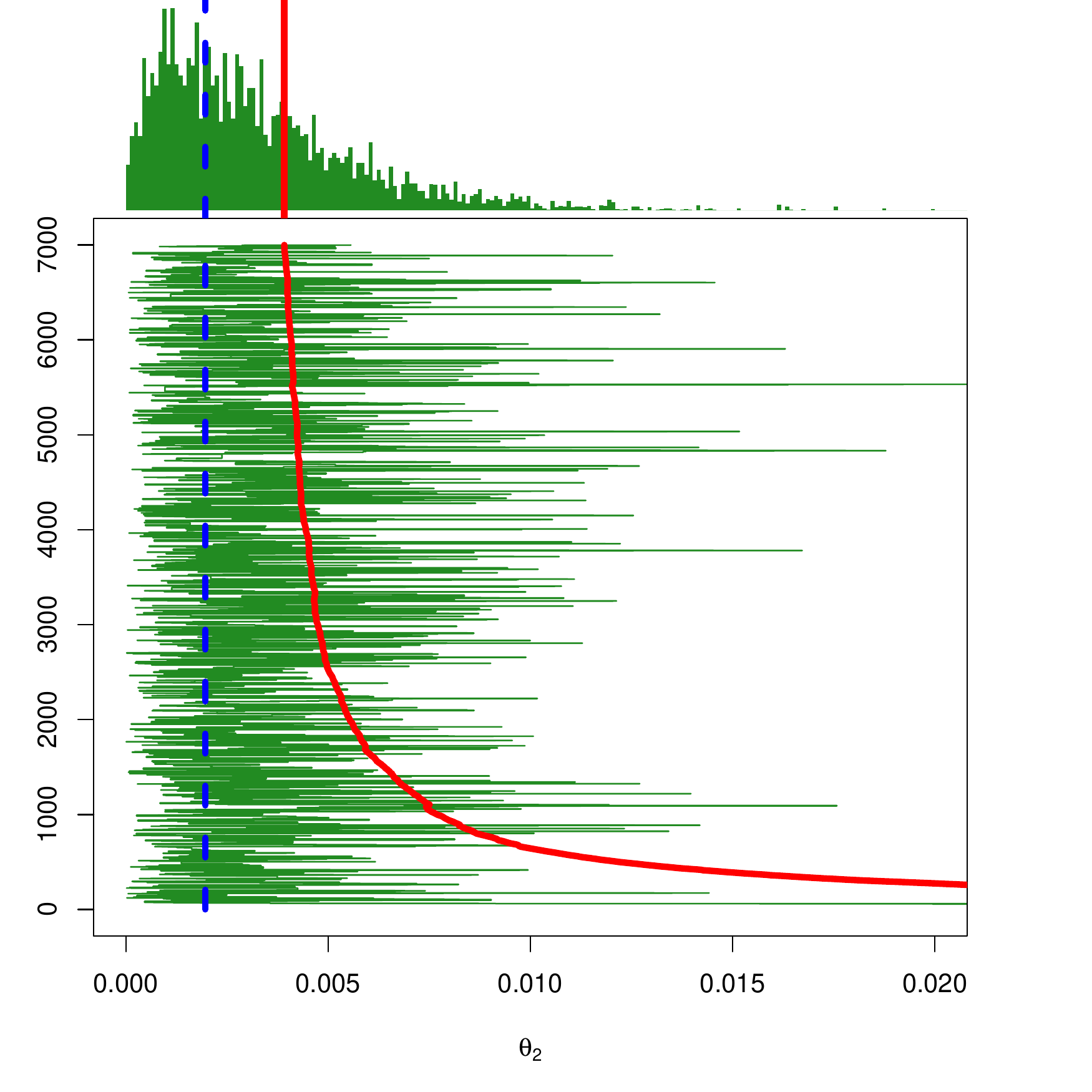}
\\
$\theta_{1}$ & $\theta_{2}$
\\[0.5em]
\rotatebox{90}{\hskip 3.5em \scriptsize MCMC iterations}
\includegraphics[width=5cm]{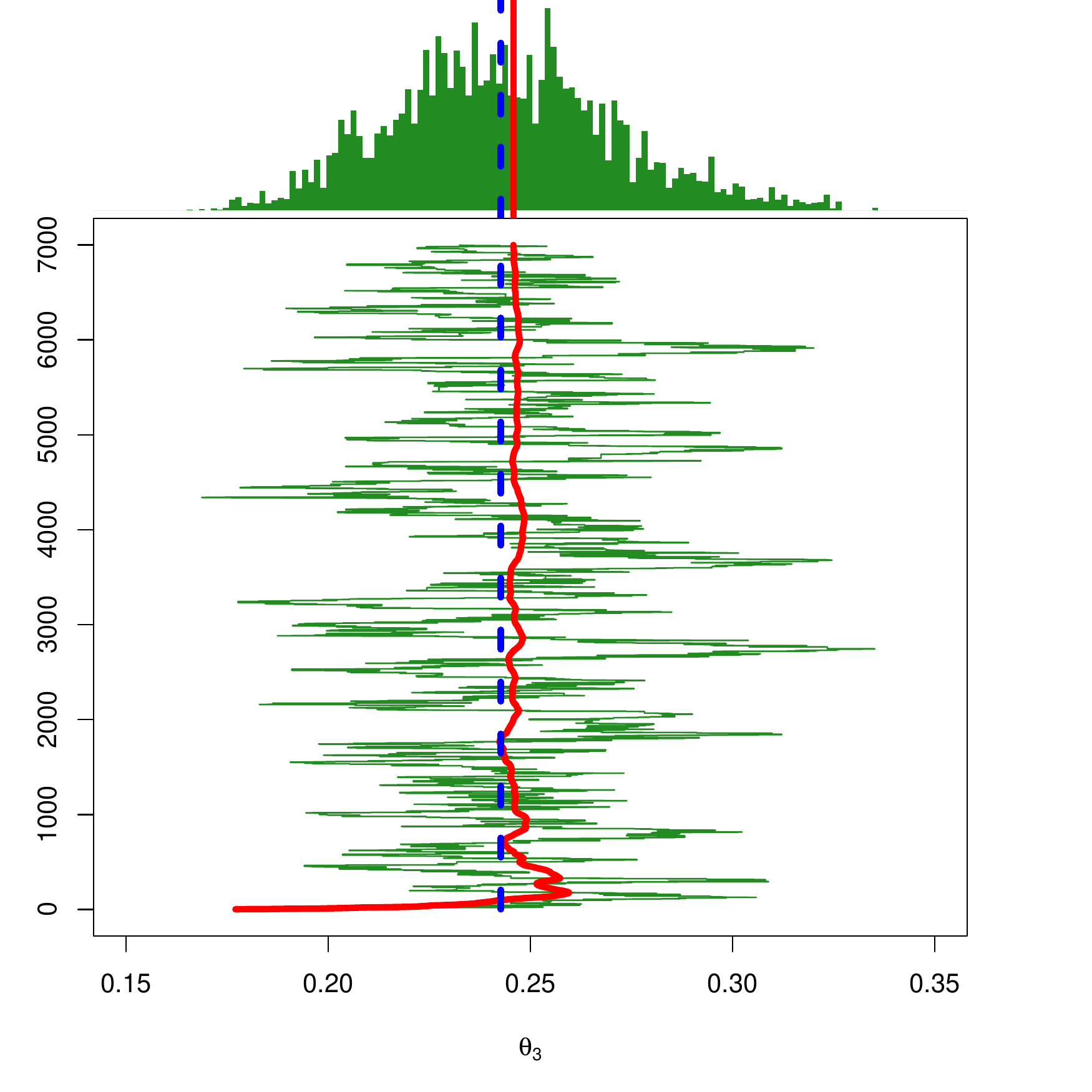}
&
\rotatebox{90}{\hskip 3.5em \scriptsize MCMC iterations}
\includegraphics[width=5cm]{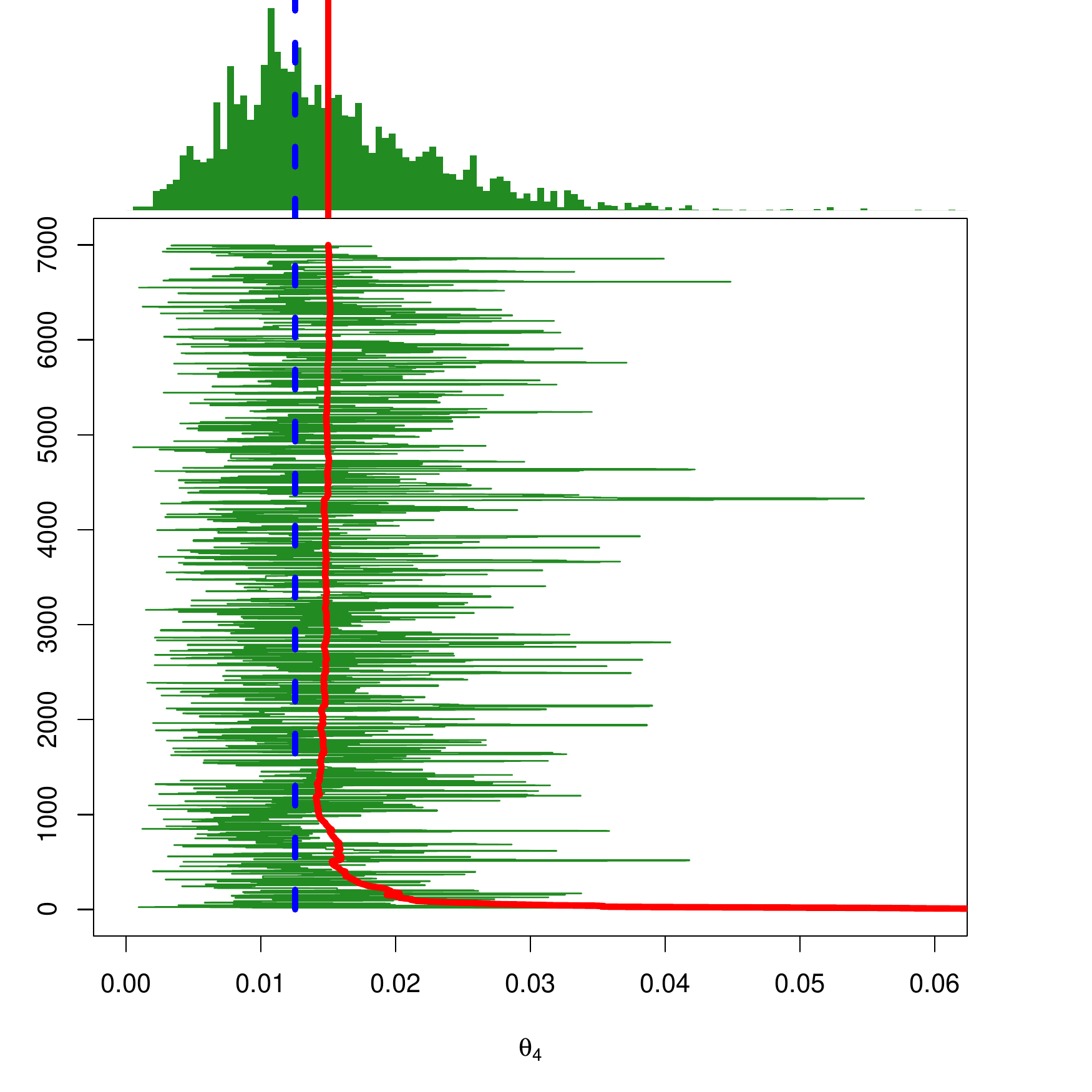}
\\
$\theta_{3}$ & $\theta_{4}$
\\[0.5em]
\multicolumn{2}{c}
{\rotatebox{90}{\hskip 3.5em \scriptsize MCMC iterations}
\includegraphics[width=5cm]{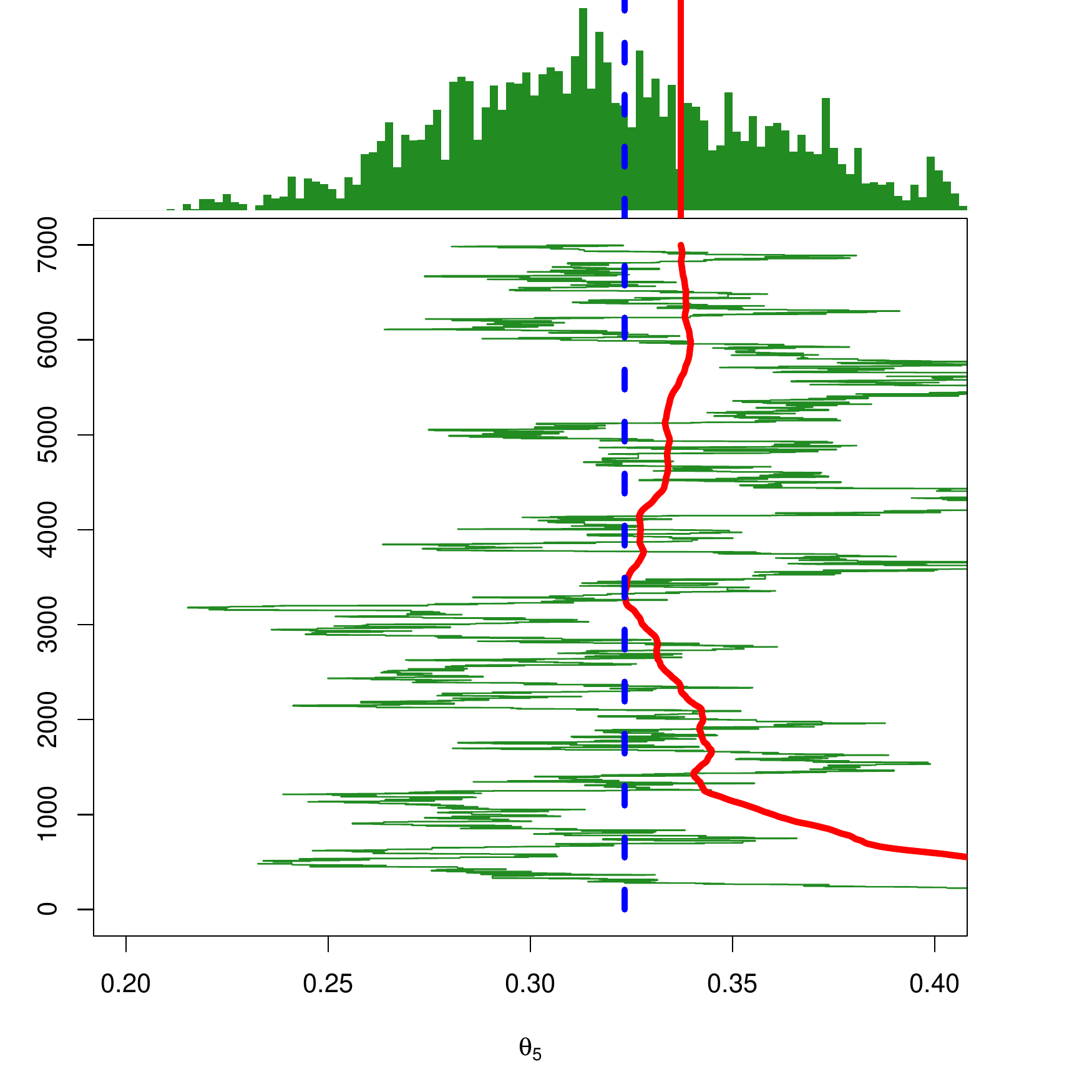}}
\\
\multicolumn{2}{c}{$\theta_{5}$}
\end{tabular}
\caption{\it MLEs and Bayes estimators for the real data set: posterior empirical distributions given by the MCMC iterations (---/green) and the associated mean (- - -/red) and the maximum likelihood estimates ($\cdots$/blue).}
\label{real_data}
\end{figure}

For the real data, the result of both approaches are slightly different (see Figure  \ref{real_data} and the Table \ref{table.estimates}). 
This calls into question the considered model. We develop this point in the next section.

%------------------
\begin{table}
\begin{center}
\small
\begin{tabular}{c|c||c|c|c|c|}
\multicolumn{2}{c}{}&\multicolumn{4}{c}{\it Bayesian estimates}\\
\multicolumn{2}{c}{}&\multicolumn{4}{c}{\it to}\\
\cline{2-6}
\multirow{6}{*}{\begin{sideways}\it from\end{sideways}}&  & $F$ & $C$ & $J$ & $B$ \\ 
\cline{2-6} \cline{2-6}
  & $F$ 
  & 0.9121 \cellcolor[gray]{0.8}  
  & 0.0842 
  & 0.0037 
  & 0 \cellcolor[gray]{0.4} \\ \cline{2-6}
  & $C$ 
  & 0    \cellcolor[gray]{0.4} 
  & 0.7417  \cellcolor[gray]{0.8} 
  & 0.2433 
  & 0.0150 \\ \cline{2-6}
  & $J$ 
  & 0      \cellcolor[gray]{0.4} 
  & 0.3273 
  & 0.6727    \cellcolor[gray]{0.8} 
  & 0  \cellcolor[gray]{0.4}\\ \cline{2-6}
  & $B$ 
  & 0 \cellcolor[gray]{0.4} 
  & 0 \cellcolor[gray]{0.4} 
  & 0 \cellcolor[gray]{0.4} 
  & \ 1 \cellcolor[gray]{0.4} \\ \cline{2-6}
\end{tabular}\\[2em]
\begin{tabular}{c|c||c|c|c|c|}
\multicolumn{2}{c}{}&\multicolumn{4}{c}{\it Maximum likelihood estimates}\\
\multicolumn{2}{c}{}&\multicolumn{4}{c}{\it to}\\
\cline{2-6}
\multirow{6}{*}{\begin{sideways}\it from\end{sideways}}&  & $F$ & $C$ & $J$ & $B$ \\ 
\cline{2-6} \cline{2-6}
  & $F$ 
  & 0.9158 \cellcolor[gray]{0.8} 
  & 0.0823 
  & 0.0019 
  & 0 \cellcolor[gray]{0.4} \\ \cline{2-6}
  & $C$ 
  & 0       \cellcolor[gray]{0.4} 
  & 0.7449  \cellcolor[gray]{0.8} 
  & 0.2426 
  & 0.0125 \\ \cline{2-6}
  & $J$ 
  & 0 \cellcolor[gray]{0.4} 
  & 0.3233 
  & 0.6767  \cellcolor[gray]{0.8}  
  & 0 \cellcolor[gray]{0.4}\\ \cline{2-6}
  & $B$ 
  & 0 \cellcolor[gray]{0.4} 
  & 0 \cellcolor[gray]{0.4} 
  & 0 \cellcolor[gray]{0.4} 
  & \ 1 \cellcolor[gray]{0.4} \\ \cline{2-6}
\end{tabular}
\end{center}
\caption{\it Bayesian and  Maximum likelihood estimates. Dark grey cells correspond to transition probabilities that are supposed to be known; light grey cells correspond to transition probabilities that are deduced from the other ones.}
\label{table.estimates}
\end{table}
%------------------

%------------------
%------------------

%-------------------------------------------------------------------------
\subsubsection*{\it Distribution of the time to reach $B$}
%-------------------------------------------------------------------------

Given the two estimations of the transition matrix, we would like to address the two following questions. First, what is the distribution law of the first time to reach the absorbing state~$B$~? Second, as $B$ is absorbing and all other states are transient, the limit distribution of the Markov chain is $\delta_{B}$, but before this state $b$ is reached what is the ``limit'' distribution of $X_{n}$ on the other states~? This distribution is called the quasi-stationary distribution of the process $X_{n}$ and we will compute it.

\begin{figure}
\centering
\includegraphics[width=5cm]{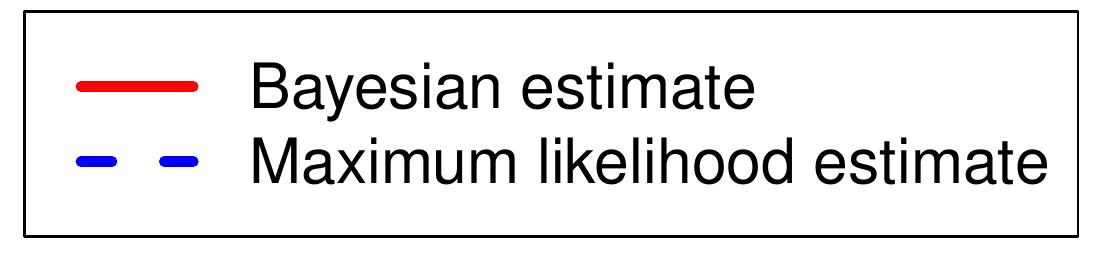}\\[0.8em]
\includegraphics[width=12cm]{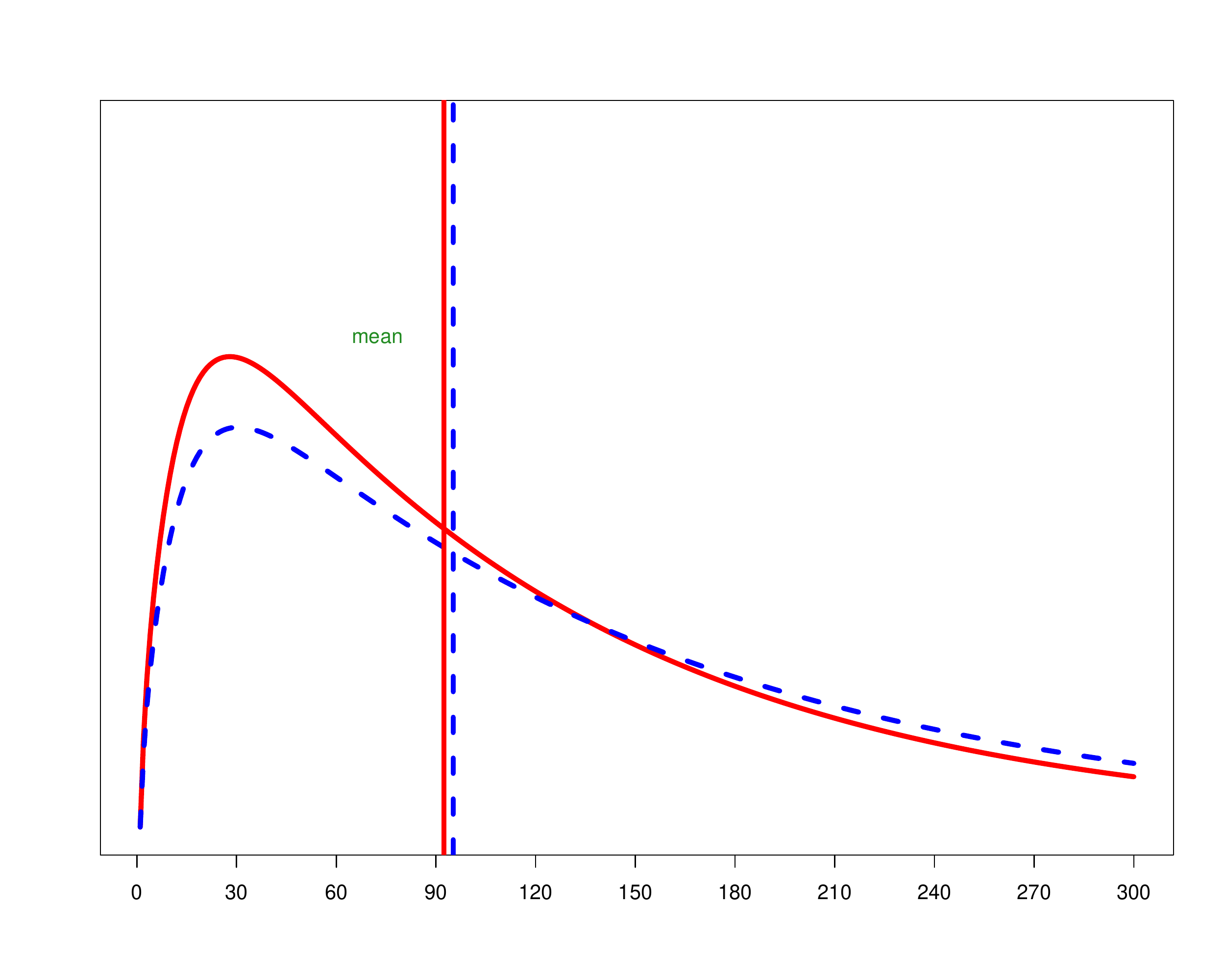}
\caption{\it Distribution of the time to reach the state ``B''. The mean time is 92 years for the Bayesian estimate and 96 years for the MLE.}
\label{fig.timeToB}
\end{figure}

To answer the first question we use the result of Appendix B: the distribution law of the first time $\tau_{FB}$ to reach $B$ starting from $F$:
\begin{align*}
  \P(\tau_{FB}=n|X_{0}=F)
  =
  \P(X_{n}=B,\, X_{m}\neq B,\, m=1,\dots,n-1|X_{0}=F)
\end{align*}
is given by recurrence formula \eqref{eq_ff} and  plotted in Figure \ref{fig.timeToB} for both the Bayesian and the maximum likelihood estimates. The mean time is 92 years for the Bayesian estimate and 96 years for the MLE.

%------------------
%------------------

%-------------------------------------------------------------------------
\subsubsection*{\it Limit distribution before reaching $B$ (quasi-stationary distribution)}
%-------------------------------------------------------------------------

The answer to the second question is given by the so-called quasi-stationary distribution, see Appendix C. From \eqref{eq.quasi.stationary} we can compute the quasi-stationary distribution $\tilde\mu=(\tilde\mu(F),\tilde\mu(C),$ $\tilde\mu(J))$ associated with the estimators of $Q$ with the maximum likelihood and the Bayesian approaches. The results are:
\begin{center}
\begin{tabular}{r|c|c|c|}
\cline{2-4}
   & $F$ & $C$ &$J$\\
\cline{2-4}
{\it Bayesian estimator} &\hskip1em0\hskip1em\mbox{}& 0.5659 &  0.4341\\
\cline{2-4}
{\it Maximum likelihood estimator} & 0 & 0.5672& 0.4328\\
\cline{2-4}
\end{tabular}
\end{center}
Hence, conditionally the fact that the process does not reach $B$, and as soon at it leaves the state $F$, it will spend $57\%$ of its time in the $C$ state and $43\%$ of its time in the $J$ state.

%%%%%%%%%%%%%%%%%%%%%%%%%%%%%%%%%%%%%%%%%%%%%%%%%%%%%%%%%%%%%%%%%%%%%%%%%%%%
%%%%%%%%%%%%%%%%%%%%%%%%%%%%%%%%%%%%%%%%%%%%%%%%%%%%%%%%%%%%%%%%%%%%%%%%%%%%
\section{Model evaluation}
\label{sec.evalution}
%%%%%%%%%%%%%%%%%%%%%%%%%%%%%%%%%%%%%%%%%%%%%%%%%%%%%%%%%%%%%%%%%%%%%%%%%%%%
%%%%%%%%%%%%%%%%%%%%%%%%%%%%%%%%%%%%%%%%%%%%%%%%%%%%%%%%%%%%%%%%%%%%%%%%%%%%

In this section we test the fit between the data and the model. From the data set of Figure \ref{fig.data} it is clear that the holding time in the state ``forest'' does not seem to correspond to that of a Markov chain. The holding time $S(e)$ of a given state $e\in E$, also called its sojourn time, is the number of consecutive time periods the Markov chain $X_{n}$ remains in this state:
\begin{align*}
  S(e) \eqdef \inf\{n\in \N;\; X_n \neq  e\}
\end{align*}
conditionally on $X_{0}=e$. The distribution law of $S(e)$ is given by:
\begin{align*}
  \P_e(S(e) = n) 
  &= \P_e(X_0=\cdots = X_n=e, X_{n+1}\neq e)
  \\
  &= \sum_{e'\neq e} \P_e (X_0=\cdots = X_n = e, X_{n+1}=e')
  \\
  &=\sum_{e'\neq e} \underbrace{Q(e,e)\cdots Q(e,e)}_{\textrm{\tiny $n$ times}}\,Q(e,e')
  = (Q(e,e))^n \, (1-Q(e,e))
\end{align*}
for $n\geq 1$ and 0 for $n=0$, that is a geometric distribution of parameter $Q(e,e)=\P(X_{n+1}=e|X_{n}=e)$. Note that $\E_e S(e)=1/(1-Q(e,e))$ and $\var_{e}(R_{e})=Q(e,e)/(1-Q(e,e))^2$.

%%%%%%%%%%%%%%%%%%%%%%%%%%%%%%%%%%%%%%%%%%%%%%%%%%%%%%%%%%%%%%%%%%%%%%%%%%%%
\subsection{Goodness-of-fit test}
%%%%%%%%%%%%%%%%%%%%%%%%%%%%%%%%%%%%%%%%%%%%%%%%%%%%%%%%%%%%%%%%%%%%%%%%%%%%

In order to test if the distribution of the holding time $S(e)$ on each state $e\in E$ of the data set $(e_n^p)_{n=0:21}^{p=1:43}$ is geometric, we use a bootstrap technique for goodness-of-fit on empirical distribution function proposed in \cite{henze1996a}.

\nocite{genest2008a,stute1993a}

Considering a sample $S_1,\dots,S_k$ of size $k$ from  a discrete cumulative distribution function $(F(n))_{n\in \N}$, we aim to test the following hypothesis:
\begin{equation}
  H_0\,:\, F \in \F \eqdef \{F_p : p\in \Theta\}\,.
\label{eq14}
\end{equation}
In our case, $F_p$ is a geometric cumulative distribution function (CDF) with parameter $p\;\in [0\;1]$. Classically, we consider an estimator:
\begin{equation*}
   \hat{p}=T(S_{1:k})
\end{equation*} 
of $p$ and we compute the distance between the theoretical CDF  $F_{\hat p}$ and the empirical CDF:
\begin{equation*}
  \hat{F}_{S_{1:k}} (n) 
  \eqdef 
  \frac{1}{k} \sum_{\ell=1}^k \indic_{\{S_\ell\leq n\}}\,.
\end{equation*}
We use the Kolmogorov-Smirnov distance defined by:
\begin{equation}
\label{eq.Kolmogorov-Smirnov}
  K^*=K(S_{1:k})
  \eqdef 
  \sup_{n\in \N} \, \sqrt{k}\,|\hat{F}_{S_{1:k}} (n)  - F_{T(S_{1:k})} (n)|
\end{equation}
To establish whether $K^*$  is significantly different from 0 or not, we simulate $M$ samples of size~$k$:
\[
   S^m_1,\dots, S^m_k \simiid F_{\hat p}, \qquad m=1\cdots M
\]
and we let:
\[
  K^m \eqdef K(S_{1:k}^m),
\]
where $K$ is the function defined in \eqref{eq.Kolmogorov-Smirnov}.

The $p$-value associated to that test is:
\[
   \rho \eqdef \frac{1}{M}\sum_{m=1}^M \indic_{\{K^m \geq K^*\}}
   \,.
\]
If $\rho$ is less than a given threshold $\alpha$, corresponding to the probability
chance of rejecting the null hypothesis $H_0$ when it is true à tort, then $H_0$ is rejected.
 
%%%%%%%%%%%%%%%%%%%%%%%%%%%%%%%%%%%%%%%%%%%%%%%%%%%%%%%%%%%%%%%%%%%%%%%%%%%%
\subsection{Holding time goodness-of-fit test}
%%%%%%%%%%%%%%%%%%%%%%%%%%%%%%%%%%%%%%%%%%%%%%%%%%%%%%%%%%%%%%%%%%%%%%%%%%%%

The $B$ state (perennial crop) is absorbing so that its holding time is infinite. Moreover, states that appear at the end of the series of Figure \ref{fig.data} are not treated (they are considered as censored data). Then the holding time values on each state $F$, $C$, $J$ in the data set are given in Table \ref{table.holding.time.values}.

%---------------------------------------------
\begin{table}
\begin{center}
\begin{tabular}{|llccccccc|}
\hline
$F$ (forest) & Holding time values& 1 & 3 & 11 & 13 & 14 & 15 & 16\\
& Number of occurrences & 1 & 9 & 2 & 1 & 6 & 21 & 3\\
\hline
\hline
$C$ (annual crop) & Holding time values& 1 & 2 & 3 & 4 & 5 & 6 &\\
& Number of occurrences & 11 & 17 & 12 & 5 & 9 & 7 &\\
\hline
\hline
$J$ (fallow) & Holding time  values& 1 & 2 & 3 & 4 & 6 & 8 & 11 \\
&  Number of occurrences & 16 & 12 & 7 & 4 & 2 & 1 & 1  \\
\hline
\end{tabular}
\end{center}
\caption{\it Holding time values (year) on each state $F$, $C$, $J$ in the data set.}
\label{table.holding.time.values}
\end{table}
%---------------------------------------------

In order to test the hypothesis $H_{0}$ we use the MLE for the parameter $p$ of the geometric PDF:
\begin{align}
\label{eq.emv.p}
 \hat{p} \eqdef \frac{1}{1 + \frac{1}{k}\sum_{\ell=1}^k S_\ell } \,.
\end{align}
Indeed, the likelihood function is:
\[
   L(p)
   =
   (1-p)^{S_1}\,p \cdots (1-p)^{S_k}\,p
   = 
   (1-p)^{\sum_{\ell=1}^{k}S_\ell}\,p^k
\]
and $L'(p)=0$ leads to
$  n-p \,(\sum_{\ell=1}^n S_\ell + n)=0$ and \eqref{eq.emv.p}.

The complete test procedure is given by Algorithm \ref{alg2}.

%------------------------
\begin{algorithm}
\begin{center}
\begin{minipage}{11cm}
\hrulefill\\[-1em]
\mbox{}
\begin{algorithmic}
\STATE  $\hat p  \ot  T(S_{1:k})$
\STATE  $\bar n\ot \sup(S_1,\ldots,S_k)$
\STATE $K^* \ot \, 
   \sup\{\sqrt{k}|F^{S_{1:k}}(n)-F^{T(S_{1:k})}(n)|
     \,,\, 0\leq n\leq \bar n\}$
\FOR  {$m=1,2,\ldots,M$}
  	\STATE  $S^m_1,\dots, S^m_k \simiid F_{\hat p}$
  	\STATE  $\bar n\ot \sup(S_1^m,\ldots,S_k^m)$
  	\STATE $K^m \ot \, 
    	 \sup\{\sqrt{k}|F^{S^m_{1:k}}(n)-F^{T(S^m_{1:k})}(n)|
    		 \,,\, 0\leq n\leq \bar n\}$
\STATE $\rho \ot  \frac{1}{M}\sum_{m=1}^M \indic_{\{K^m\geq K^*\}}$
\ENDFOR

\IF{$\rho\leq\alpha$}
   \STATE accept $H_0$,
\ELSE
  \STATE reject $H_0$. 
\ENDIF
\end{algorithmic}
\hrulefill
\end{minipage}
\end{center}
\caption{\it Parametric bootstrap for goodness-of-fit with the geometric distribution of parameter $p$.}
\label{alg2}
\end{algorithm}
%------------------------

%%%%%%%%%%%%%%%%%%%%%%%%%%%%%%%%%%%%%%%%%%%%%%%%%%%%%%%%%%%%%%%%%%%%%%%%%%%%
\subsection{Results}
%%%%%%%%%%%%%%%%%%%%%%%%%%%%%%%%%%%%%%%%%%%%%%%%%%%%%%%%%%%%%%%%%%%%%%%%%%%%

%------------------
\begin{figure}
 \centering
\begin{tabular}{ccc}
\includegraphics[width=4cm]{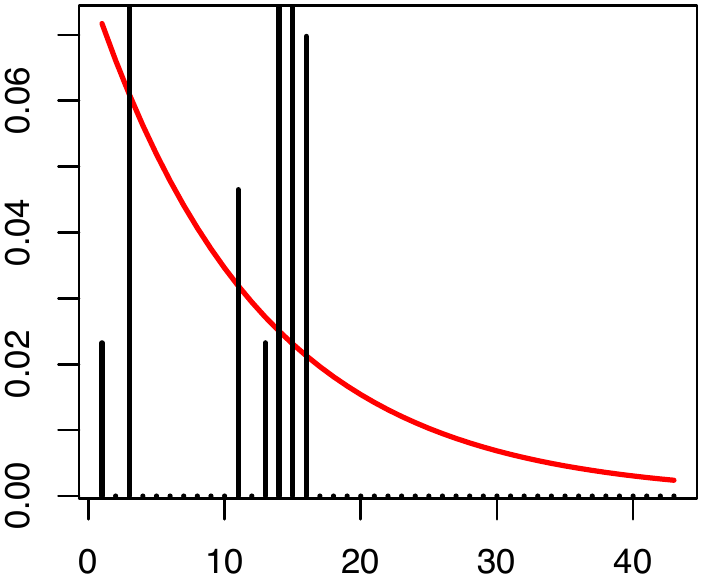}
&
\includegraphics[width=4cm]{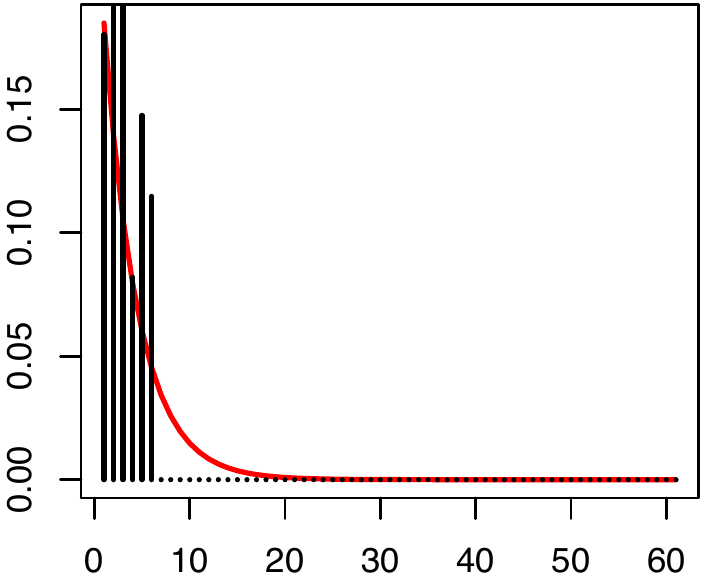}
&
\includegraphics[width=4cm]{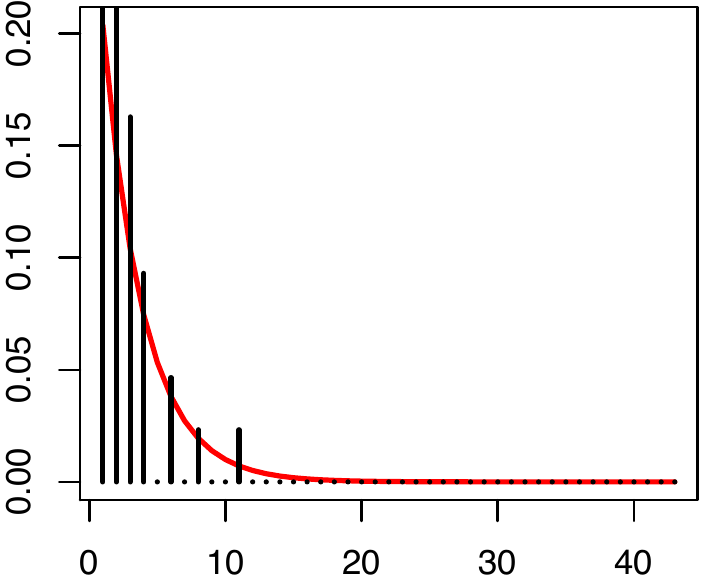}
\\
\small Forest $F$ & \small Annual crop $C$ & \small Fallow $J$
\end{tabular}
\caption{\it Empirical PDFs of holding time of states $(F,C,J)$ associated with the data set of Figure \ref{fig.data}; and the geometric PDF, represented as continuous read lines, corresponding to the parameter $p$ estimated by \eqref{eq.emv.p}.}
\label{fig.Tps_sej}
\end{figure}
%------------------

\begin{figure}
 \centering
\begin{tabular}{ccc}
\includegraphics[width=4cm]{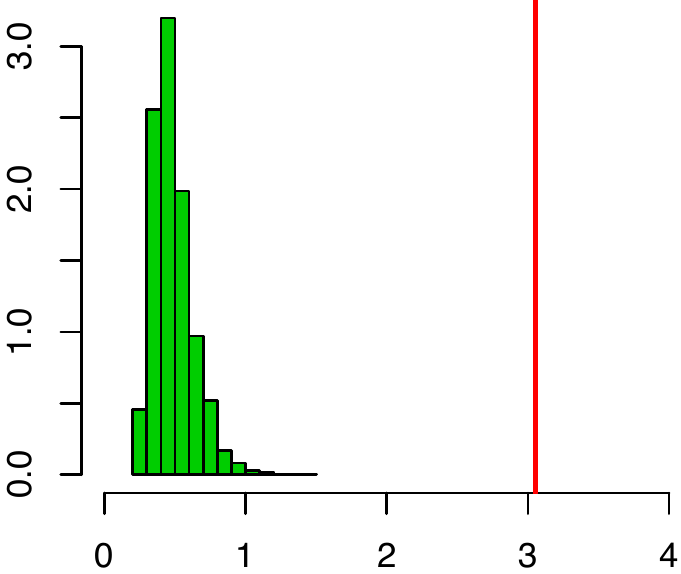}
&
\includegraphics[width=4cm]{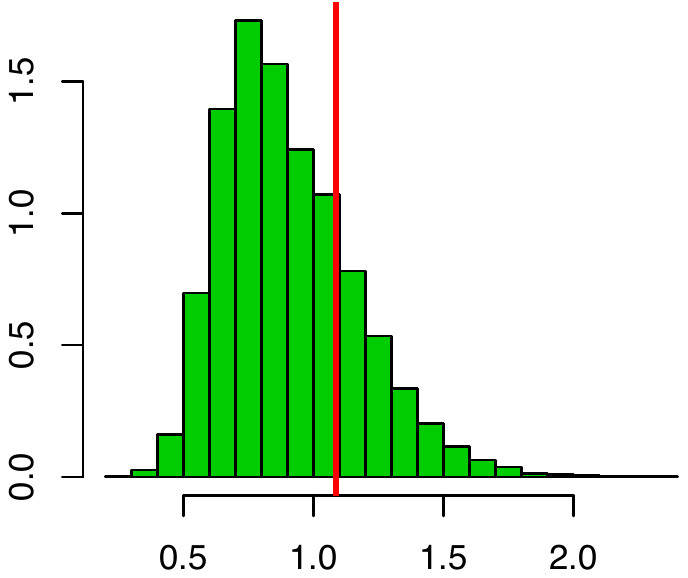}&
\includegraphics[width=4cm]{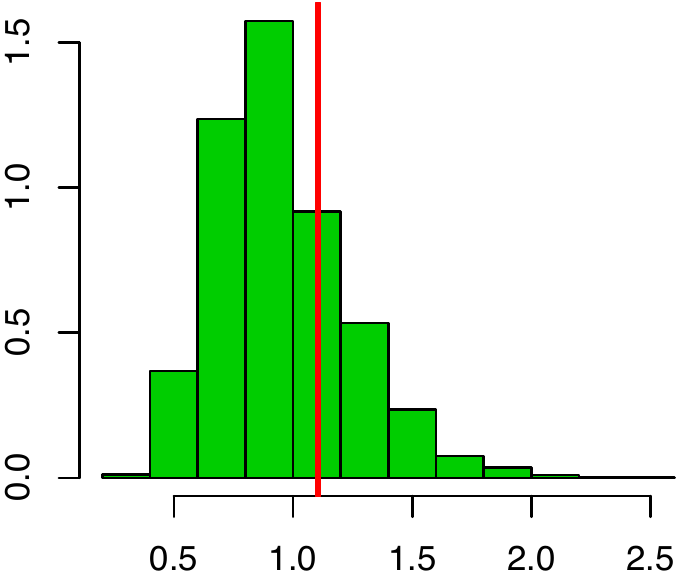}
\\
\small Forest $F$ & \small Annual crop $C$ & \small Fallow $J$
\end{tabular}
\caption{\it Empirical PDF of $K^m$ (sample from $K^*$) for the three states $(F,C,J)$  et  $K^*$ value (vertical line).}
\label{fig.boot}
\end{figure}

In Figure \ref{fig.Tps_sej} we plotted  the empirical PDFs of holding time of states ``Forest'', ``Annual crop'' and ``Fallow'', associated with the data set of Figure \ref{fig.data}, and the geometric PDF corresponding to the parameter $p$ estimated by \eqref{eq.emv.p}. We see that in the case of the ``Forest'' state the matching is questionable. In Figure \ref{fig.boot} we plotted the empirical PDF of $K^m$ for the three states. The value of $K^*$ and of the associated $p$-values are:
\begin{center}
\begin{tabular}{|l||c|c|c|}
\hline
   & ``Forest'' $F$ &  ``Annual crop'' $C$ & ``Fallow'' $J$\\
\hline
 $K^*$  & 3.051 & 1.086  & 1.104 \\%3.051060 1.085909 1.104389
\hline
$p$-value &  0  &  0.224 &  0.255\\
\hline
\end{tabular}
\end{center} 

%------------------
%------------------

In conclusion, the geometric distribution hypothesis is strongly rejected for the state ``forest''. The $p$-value for this state is null. This is understandable as this state does not really correspond to a ``dynamic state''.

%%%%%%%%%%%%%%%%%%%%%%%%%%%%%%%%%%%%%%%%%%%%%%%%%%%%%%%%%%%%%%%%%%%%%%%%%%%%
%%%%%%%%%%%%%%%%%%%%%%%%%%%%%%%%%%%%%%%%%%%%%%%%%%%%%%%%%%%%%%%%%%%%%%%%%%%%
\section{Conclusion and perspectives}
\label{sec.conlusion}
%%%%%%%%%%%%%%%%%%%%%%%%%%%%%%%%%%%%%%%%%%%%%%%%%%%%%%%%%%%%%%%%%%%%%%%%%%%%
%%%%%%%%%%%%%%%%%%%%%%%%%%%%%%%%%%%%%%%%%%%%%%%%%%%%%%%%%%%%%%%%%%%%%%%%%%%%

We proposed a Markovian model of land use dynamics for parcels near
the forest corridor of Ranomafana and Andringitra national parks in
Madagascar. We supposed first the dynamic uses of the parcels are 
independent and identically distributed; second that the dynamics is
Markovian with four states. The transition matrix depends on five
unknown parameters. We considered the MLE and the Bayes estimate with
Jeffrey prior. In this last case, the estimator is computed with a
MCMC procedure. The Bayes estimator performs slightly better than the
MLE. On the real data set, the two estimators give rather similar
results.

We assessed the adequacy of the model to real data. We focused on the
holding times: we tested if the empirical holding times correspond to
a geometric distribution.  
We used a parametric bootstrap goodness-of-fit on empirical
distribution. Clearly the geometric distribution hypothesis is
violated in the case of the ``Forest'' state. 

The ``Forest'' state therefore requires a special treatment. In a near
future we are now developing a semi-Markov model where the sojourn
time on the state $F$ will better match the data set and so will not
be geometric.

The long time behavior of the inferred model is dubious as the present
data set is relatively limited in time (22 years).  This data set
implies a relatively short time scale where some rare transitions,
like the forest regeneration, are not observed. Note that the
Bayesian approach has an advantage over the likelihood approach in
that it allows to incorporate prior knowledge about these rare and
unobserved transitions. The likelihood approach will set their
probabilities zero while the Bayesian approach will incorporate a
priori knowledge and assign them positives probabilities. A new
database is currently being developed by the IRD. It will be for a
longer period of time and a greater number of parcels, it will also
allow to consider a more detailed state space comprising more than
four states. In a longer time scale, it is reasonable to suppose that $F$ and $B$
have long sojourn time distributions, the one associated to $F$ being
longer than the one associated to $B$. Also $B$ will not be absorbing
anymore as well as the forest regeneration will be possible, i.e. the
transition from $J$ to $F$ will be possible. The associated model will present
multi-scale properties, namely slow and fast components in the dynamics, that will be of interest.

Part of the complexity of these agro-ecological temporal data comes
from the fact that some transitions are ``natural'' while others come
from human decisions (annual cropping, crop abandonment, planting perennial crops,
etc.). It should also be interesting to study the dynamics of parcels
conditionally on the dynamics of the neighbor parcels. This model
could be more realistic but requires first studying the farmers'
practices in order to limit the number of unknown parameters in the
model.

\newpage
\appendix

%%%%%%%%%%%%%%%%%%%%%%%%%%%%%%%%%%%%%%%%%%%%%%%%%%%%%%%%%%%%%%%%%%%%%%
%%%%%%%%%%%%%%%%%%%%%%%%%%%%%%%%%%%%%%%%%%%%%%%%%%%%%%%%%%%%%%%%%%%%%%
\section*{Appendices}
\addcontentsline{toc}{section}{Appendices}
%%%%%%%%%%%%%%%%%%%%%%%%%%%%%%%%%%%%%%%%%%%%%%%%%%%%%%%%%%%%%%%%%%%%%%
%%%%%%%%%%%%%%%%%%%%%%%%%%%%%%%%%%%%%%%%%%%%%%%%%%%%%%%%%%%%%%%%%%%%%%

%%%%%%%%%%%%%%%%%%%%%%%%%%%%%%%%%%%%%%%%%%%%%%%%%%%%%%%%%%%%%%%%%%%%%%
\addcontentsline{toc}{subsection}{A. Explicit Bayes estimators for the two state case}
\subsection*{A. Explicit Bayes estimators for the two state case}
%%%%%%%%%%%%%%%%%%%%%%%%%%%%%%%%%%%%%%%%%%%%%%%%%%%%%%%%%%%%%%%%%%%%%%

Let  $(X_n)_{0\leq n\leq N}$  be a Markov chain with two states $\{0,1\}$ and transition matrix 
\[ Q \eqdef 
       \left(\begin{smallmatrix}
         p    & 1-p  \\
         1-q  & q
       \end{smallmatrix}\right)\,.
\]
We suppose that the initial law is the invariant distribution $\mu= ({q}/{p + q},{p}/{p+q})$, that is the solution of  $\mu\,Q=\mu$. The unknown parameter is $\theta = (p,q) \in [0,1]^{2}$ and the associated likelihood function is
\[
   L_N(\theta) 
   \eqdef
   \P_{\theta}(X_{0:N} = x_{0:N}) 
   = (1-p)^{n_{00}}\, p^{n_{01}} \,q^{n_{10}} \,(1-q)^{n_{11}}
\]
where
$  n_{ij}
  \eqdef
  n_{ij}(x_{0:N})
  =
  \sum^{N-1}_{n=0}\indic_{\{X_n=i\}}\,\indic_{\{X_{n+1}=j\}}
$
is the number of transition $i\to j$ in $x_{0:N}$.

We consider the following priori distributions: the uniform distribution $\pi^{\textrm{\rm\tiny U}}$ on $[0,1]^2$ and the beta distribution $\pi^{\textrm{\rm\tiny B}}$ with parameters $(a,b)$, that is
\[
  \pi^{\textrm{\rm\tiny B}}(\theta)
  =
  \frac{1}{\beta(a,b)}\;\theta^{a-1}\,(1-\theta)^{b-1}
\]
where $\beta (a,b)$ is the beta function:
\[
 \beta (a,b)
 \eqdef
 \int_0^1 x^{a-1}\, (1-x)^{b-1}\,\rmd x
 =
 \frac{ \Gamma (a)\,\Gamma (b)}{\Gamma (a+b)}
\]
with
$ \Gamma(z) \eqdef \int_0^{+\infty} t^{z-1}\,e^{-t}\,\rmd t
$. Here we will choose $a=b=1/2$, note that $\Gamma(1)=1$ and $\Gamma(\frac{1}{2})=\sqrt{\pi}$.
For these two priors we can explicitly compute the posterior distriution and the associated Bayes estimators. Indeed the posterior distribution $\tilde{\pi}$ is given by the Bayes formula:
$ \tilde{\pi}(\theta)  \propto L_N(\theta)\,\pi(\theta)$, that is:
\begin{align*}
  \tilde{\pi}^{\textrm{\rm\tiny U}} (\theta) 
  & \propto 
  L_N (\theta)
  =
  (1 - p)^{n_{00}}\, p^{n_{01}}\,q^{n_{10}}\,(1 - q)^{n_{11}}\,,  
\\
  \tilde{\pi}^{\textrm{\rm\tiny B}} (\theta) 
  & \propto 
  L_N (\theta)\,\pi^{\textrm{\rm\tiny B}}(\theta) 
  = 
  (1 - p)^{n_{00} - \frac{1}{2}}\, p^{n_{01} - \frac{1}{2}}
    \,q^{n_{10} - \frac{1}{2}}\,(1 - q)^{n_{11} - \frac{1}{2}}
\end{align*}
and the corresponding  Bayes estimator are:
\begin{align*}
   \tilde{\theta}^{\textrm{\rm\tiny U}} 
   &= 
   \int_{[0,1]^2} \theta\; \tilde{\pi}^{\textrm{\rm\tiny   U}}(\theta) \,\rmd \theta
   \,,
&
   \tilde{\theta}^{\textrm{\rm\tiny B}} 
   &= 
   \int_{[0,1]^2} \theta \;\tilde{\pi}^{\textrm{\rm\tiny B}}(\theta) \,\rmd\theta\,.
\end{align*}
We can easily check that the estimators of $p$ and $q$ for the uniform prior:
\begin{align*}
  \tilde{p}^{\textrm{\rm\tiny U}} 
  &= 
  \frac{n_{01} + 1}{n_{01} + n_{00} + 2}\,, 
&
  \tilde{q}^{\textrm{\rm\tiny U}} 
  &= 
  \frac{n_{10} + 1}{n_{10} + n_{11} + 2}
\end{align*}
and for the beta prior:
\begin{align*}
  \tilde{p}^{\textrm{\rm\tiny B}} 
  &= 
  \frac{n_{01} + \frac{1}{2}}{n_{01} + n_{00} + 1}\,,
&
  \tilde{q}^{\textrm{\rm\tiny B}} 
  &= 
  \frac{n_{10} + \frac{1}{2}}{n_{10} + n_{11} + 1}\,.
\end{align*}
Note that in this case the MLE estimators are:
\begin{align*}
  \hat{p}^{\textrm{\rm\tiny MLE}} 
  & = 
  \frac{n_{01}}{n_{00} + n_{01}}
  \,,
&
 \hat{q}^{\textrm{\rm\tiny EMV}} 
 & = 
 \frac{n_{10}}{n_{11} + n_{10}}
 \,.
\end{align*}

%%%%%%%%%%%%%%%%%%%%%%%%%%%%%%%%%%%%%%%%%%%%%%%%%%%%%%%%%%%%%%%%%%%%%%
\addcontentsline{toc}{subsection}{B. Distribution law of the time to reach a given state}
\subsection*{B. Distribution law of the time to reach a given state}
%%%%%%%%%%%%%%%%%%%%%%%%%%%%%%%%%%%%%%%%%%%%%%%%%%%%%%%%%%%%%%%%%%%%%%

Let $X_{n}$ be an homogeneous Markov chain with finite state space $E$ and transition matrix $Q$. We aim to get an explicit expression of the distribution law  $f_{ee'}^{(n)} \eqdef \P (\tau_{ee'}=n)$, $n\in \N$, of the first time $\tau_{ee'}$ to reach state $e'$ after leaving state $e$ defined as:
\[
  \tau_{ee'} 
  \eqdef
  \inf\{n\geq 1\,:\ X_n=e'|X_0=e\}\,.
\]
For $n>1$:
\begin{align*}
  Q^{(n)}(e,e')
  &=
  \P(X_{n}=e'|X_{0}=e)
  \\ 
  &=
  \P(X_{n}=e',\, \tau_{ee'}=1|X_{0}=e)
  +
  \P(X_{n}=e',\, \tau_{ee'}=1|X_{0}=e)
  +\cdots
  \\
  &\qquad\qquad\qquad
  \cdots
  +
  \P(X_{n}=e',\, \tau_{ee'}=n-1|X_{0}=e)  
  +
  \P(X_{n}=e',\, \tau_{ee'}=n|X_{0}=e)  
  \\ 
  &=
  f_{ee'}^{(1)} Q^{(n-1)}(e',e') + f^{(2)}_{ee'}\, Q^{(n-2)}(e',e') 
  + \cdots 
  + f^{(n-1)}_{ee'} Q^{(1)}(e',e')
  + f_{ee'}^{(n)} 
\end{align*}
hence $f^{(n)}_{ee'}$ could be computed recursively according to
\begin{align}
\label{eq_ff}
  f^{(n)}_{ee'}
  &= Q^{(n)}(e,e') - \sum_{k=1}^{n-1}f^{(k)}_{ee'} Q^{(n-k)}(e',e') 
\end{align}
with $f_{ee'}^{(1)}= Q^{(1)}(e,e')$.

%%%%%%%%%%%%%%%%%%%%%%%%%%%%%%%%%%%%%%%%%%%%%%%%%%%%%%%%%%%%%%%%%%%%%%
\addcontentsline{toc}{subsection}{C. Quasi-stationary distribution}
\subsection*{C. Quasi-stationary distribution}
%%%%%%%%%%%%%%%%%%%%%%%%%%%%%%%%%%%%%%%%%%%%%%%%%%%%%%%%%%%%%%%%%%%%%%

We consider the probability to be in $e\in\{F,C,J\}$ before reaching $B$ and starting from $F$:
\begin{align*}
\mu_{n}(e)
  &\eqdef 
  \P(X_n=e|X_m\neq B,\,m=1,\dots,n-1,\, X_0=F)
\\
  &
  =
  \P(X_n=e|X_{n-1}\neq B,\,X_0=F) 
\\
  &
  =
  \frac{\P(X_n=e,\,X_{n-1}\neq B|X_0=F)}
       {\P(X_{n-1}\neq B|X_0=F)},
\\
  &
  =
  \frac{\P(X_n=e|X_0=F)}
       {1-\sum_{e'\in\{F,C,J\}}\P(X_{n-1}=e|X_0=F)}\,.
\end{align*}
When
\begin{align*}
  \mu_{n}(e) \xrightarrow[n\to\infty]{}\tilde\mu(e)\,,\ e\in\{F,C,J\}
\end{align*}
the probability distribution $(\tilde\mu(e))_{e\in\{F,C,J\}}$ is called quasi-stationary probability distribution. This problem was originally solved in \cite{darroch1965a}: $\tilde\mu=[\tilde\mu(F)\ \tilde\mu(C)\ \tilde\mu(J)]$ exists and it is given by the equation
\begin{align}
\label{eq.quasi.stationary}
   \tilde\mu \,\tilde Q = \tilde\lambda\,\tilde\mu
\end{align}
with $\tilde\mu(e)\geq 0$ and $\tilde\mu(F)+\tilde\mu(C)+\tilde\mu(J)=1$, where $\tilde Q$ is the $3\times 3$ submatrix defined by
\begin{align*}
Q = \left(
         \begin{array}{c|c}
             \tilde Q & \tilde q \\ \hline
             0 & 1
         \end{array}
      \right ),
\end{align*}
and $\tilde\lambda$ is the spectral radius of $\tilde Q$.

%%%%%%%%%%%%%%%%%%%%%%%%%%%%%%%%%%%%%%%%%%%%%%%%%%%%%%%%%%%%%%%%%%%%%%%%%%%%%%%%%%%%
%%%%%%%%%%%%%%%%%%%%%%%%%%%%%%%%%%%%%%%%%%%%%%%%%%%%%%%%%%%%%%%%%%%%%%%%%%%%%%%%%%%%
%%%%%%%%%%%%%%%%%%%%%%%%%%%%%%%%%%%%%%%%%%%%%%%%%%%%%%%%%%%%%%%%%%%%%%%%%%%%%%%%%%%%
\newpage
\nocite{kemeny1976a}
\addcontentsline{toc}{section}{References}
%\bibliography{lib/fab}
\bibliographystyle{plain}%plain

%%%%%%%%%%%%%%%%%%%%%%%%%%%%%%%%%%%%%%%%%%%%%%%%%%%%%%%%%%%%%%%%%%%%%%%%%%%%%%%%%%%%
%%%%%%%%%%%%%%%%%%%%%%%%%%%%%%%%%%%%%%%%%%%%%%%%%%%%%%%%%%%%%%%%%%%%%%%%%%%%%%%%%%%%
%%%%%%%%%%%%%%%%%%%%%%%%%%%%%%%%%%%%%%%%%%%%%%%%%%%%%%%%%%%%%%%%%%%%%%%%%%%%%%%%%%%%

%%%%%%%%%%%%%%%%%%%%%%%%%%%%%%%%%%%%%%%%%%%%%%%%%%%%%%%%%%%%%%%%%%%%%%%%%%%%%%%%%%%%
%%%%%%%%%%%%%%%%%%%%%%%%%%%%%%%%%%%%%%%%%%%%%%%%%%%%%%%%%%%%%%%%%%%%%%%%%%%%%%%%%%%%
\end{document}